\documentclass[reqno]{amsart}
\usepackage{graphicx} 
\numberwithin{equation}{section}
\usepackage{amsmath}
\usepackage{amssymb}
\usepackage{amsthm}
\usepackage{color}
\usepackage{enumerate}
\theoremstyle{definition}
\newtheorem{thm}{Theorem}[section]
\newtheorem{prop}[thm]{Proposition}

\newtheorem{rem}[thm]{Remark}
\newtheorem{lemma}[thm]{Lemma}

\newtheorem{cor}[thm]{Corollary}

\usepackage{mathrsfs}

\newcommand{\dist}{\text{dist}}
\newcommand{\hyperbolic}{\mathcal{G}}
\newcommand{\bad}{\mathcal{B}}

\title[Birkhoff spectrum for smooth interval maps]{Uniform hyperbolic
approximation and conditional variational principle of the Birkhoff spectrum for smooth interval maps}
\author{Yuya Arima}
\date{\today}

\address{Graduate School of Mathematics, Nagoya University,
Furocho, Chikusa-ku, Nagoya, 464-8602, JAPAN} 
\email{yuya.arima.c0@math.nagoya-u.ac.jp}

\subjclass[2020]{28A78, 37C05, 37D25, 37E05}
\thanks{{\it Keywords}: Multifractal analysis, Birkhoff averages, Interval maps}

\begin{document}

\begin{abstract}
In this paper, we perform a multifractal analysis of
Birkhoff averages for topologically exact $C^{1+\epsilon}$
interval maps with non-flat critical points. 
Our first main result establishes conditional
variational formulas for the Birkhoff spectrum
of a continuous function. Our second main result establishes the locally uniform
approximation of the Birkhoff spectrum on a certain
domain by the Birkhoff spectra associated with
a single sequence of uniformly expanding subsystems.
The main novelty of this paper is a new construction of induced systems. More precisely, given any sequence of invariant measures $\{\mu_n\}_{n\in\mathbb{N}}$ with positive entropy, we construct an increasing sequence of uniformly expanding induced systems $\{F_N\}_{N\in\mathbb{N}}$ with finitely many full branches over a common base interval such that each $F_N$ simultaneously captures the information carried by the measures $\mu_1,\ldots,\mu_N$. By choosing the sequence of invariant measures appropriately, we obtain a single sequence of induced systems that allows us to approximate every invariant measure with positive entropy in the weak topology while controlling its entropy and its Lyapunov exponent.
\end{abstract}

\maketitle

\section{Introduction}
Let $f:[0,1]\rightarrow [0,1]$ be a continuous map and let $\phi$ be a continuous potential. The Birkhoff average of $\phi$ at $x\in [0,1]$ is defined by
\[
\lim_{n\to\infty}\frac{1}{n}\sum_{k=0}^{n-1}\phi(f^k(x))
\]
whenever the limit exists. 
Birkhoff averages provide a way to characterize the dynamical system $f$.  For a given $f$-invariant ergodic Borel probability measure $\mu$ on $[0,1]$,
Birkhoff's ergodic theorem yields that  for $\mu$-almost every $x\in [0,1]$ the Birkhoff average of $\phi$ at $x$ converges to the space average $\int \phi d\mu$.
Thus, for $\alpha\neq \int \phi d\mu$ the set $B(\alpha)$ of points at which the Birkhoff average of $\phi$ converges to $\alpha$ is negligible with respect to $\mu$.
However, there is still a possibility that $B(\alpha)$ might be a large set from another point of view.
Moreover, if $f$ is topologically transitive,
then $B(\alpha)$ is dense in $[0,1]$ whenever
it is nonempty. To quantify the size of $B(\alpha)$, we use the Hausdorff
dimension $\dim_H(\cdot)$ with respect to the Euclidean
metric on $\mathbb{R}$ (see, for example,  \cite[Definition 15.3.3]{Urbanskinoninvertible} for the definition of the Hausdorff dimension). The study of the Birkhoff spectrum
\[
\alpha\mapsto b(\alpha):=\dim_H(B(\alpha))
\]
is a topic in multifractal analysis.   
We refer the reader to the books by Pesin \cite{Pesinbook} and Barreira \cite{Barreirabook} for an introduction to multifractal analysis.

For uniformly expanding systems, Barreira and Saussol
\cite{BarreiraSaussol} established a conditional variational
principle for the Birkhoff spectrum of a H\"older continuous
potential. This conditional variational principle was extended to the setting of one-dimensional
parabolic Markov maps by Johansson et al.
\cite{JJTOPnonuniformly}, Climenhaga \cite{Climenhagasimultaneous} and Jaerisch and Takahasi
\cite{MixedJT}.

For a large class of interval maps with critical points, the Lyapunov
spectrum and the local dimension spectrum have been studied by Iommi and Todd \cite{IommiTodd2011dimension}, and by
Gelfert, Przytycki, and Rams \cite{GelfertPrzytyckiRams}.  However, these studies rely on the fact that the spectra can be expressed as Legendre transforms of thermodynamic quantities. Consequently, it is difficult to extend their methods to the Birkhoff spectrum $b$ of general continuous potentials, since such a representation cannot be expected to hold in general. The entropy spectrum for Birkhoff averages
has been studied
(see Takens and Verbitskiy \cite{TakensVerbitskiy} and Climenhaga \cite{Climenhaga}). However, in the presence of critical points,
it is difficult to extend the arguments for
the entropy spectrum of Birkhoff averages
to the (Hausdorff dimension) Birkhoff spectrum.

For the Birkhoff spectrum of continuous
potentials, Chung \cite{Chung} established
a conditional variational principle for topologically mixing
$C^2$ interval maps under additional
hyperbolicity assumptions.
Chung and Takahasi \cite{ChungTakahasi}
established a conditional variational principle
for Benedicks-Carleson quadratic maps,
which satisfy the Collet-Eckmann condition
and a suitable recurrence condition on the
critical orbit (see also
\cite{AlejandroFernando}). 

In this paper, for topologically exact $C^{1+\epsilon}$
interval maps with non-flat critical points,
we establish conditional variational formulas for the
Birkhoff spectrum of a continuous function.

It is natural to ask whether objects associated with
non-hyperbolic dynamical systems can be approximated
by their counterparts in uniformly hyperbolic subsystems.
We investigate this question for the Birkhoff spectrum in the setting of topologically exact $C^{1+\epsilon}$
interval maps with non-flat critical points.
More precisely, we ask whether, for $\alpha\in\mathbb{R}$,
the equality
\[
b(\alpha)=
\sup\left\{
\dim_H(Y\cap B(\alpha)):
Y\subset[0,1],\ f(Y)\subset Y,\
f|_Y\text{ is uniformly expanding}
\right\}
\]
holds. In this paper, we show that there exists an increasing
sequence $\{Y_N\}_{N\in\mathbb{N}}$ of $f$-invariant
uniformly expanding subsets of $[0,1]$ such that
for every $\alpha$ in a certain domain $\mathcal{D}$,
we have
\[
\lim_{N\to\infty}\dim_H(B(\alpha)\cap Y_N)=b(\alpha).
\]
Importantly, the sequence $\{Y_N\}_{N\in\mathbb{N}}$
can be chosen independently of $\alpha\in\mathcal{D}$.
This enables us to show that the above convergence
is in fact locally uniform on $\mathcal{D}$.

The main novelty of this paper lies in a new construction of induced systems based on a version of Katok's lemma with a uniform base set. A distinctive feature of our construction is that it starts with an arbitrarily prescribed sequence of invariant measures with positive entropy and produces an increasing sequence of uniformly expanding induced systems with finitely many full branches over a common base interval. Each induced system simultaneously captures the relevant statistical and geometric information carried by the first finitely many measures in the sequence. The full-branch structure allows branches associated with different invariant measures to be freely concatenated. 

By choosing the sequence of invariant measures appropriately, we obtain a single sequence of induced systems capable of approximating the information carried by every invariant measure with positive entropy. This construction provides a framework for studying the measure-theoretic and multifractal structures of non-uniformly hyperbolic interval maps through uniformly expanding induced systems. As an application of this construction, we establish a locally uniform hyperbolic approximation of the Birkhoff spectrum. Our version of Katok's lemma with a uniform base set also plays a key role in establishing the conditional variational formulas for the Birkhoff spectrum.

\subsection{Main results}\label{sec Main results}
Let $I:=[0,1]$ and let $f:I\rightarrow I$ be a continuously differentiable map.
A critical point of $f$ is a point at which the derivative of $f$ vanishes. We denote by $\text{Crit}(f)$ the set of critical points of $f$. A critical point $c\in \text{Crit}(f)$ is said to be non-flat if there exist $l>1$ and diffeomorphisms $\varphi$ and $\psi$ on $\mathbb{R}$ such that $\varphi(c)=\psi(f(c))=0$ and for every $x$ in a neighborhood of $c$ we have 
\[
|\psi\circ f(x)|=|\varphi(x)|^l.
\]
Since $f$ is a continuously differentiable map, if $f$ has only non-flat critical points then 
\[
\#\text{Crit}(f)<\infty.
\]
The map $f$ is said to be topologically exact if for every nonempty open subset $U$ of $I$ there exists $n\in \mathbb{N}$ such that $f^n(U)=I$. 

We denote by $M(f)$ the set of $f$-invariant Borel probability measures on $I$. For $\mu\in M(f)$ we denote by $h(\mu)$ the measure-theoretic entropy defined in \cite{walters2000introduction}. Moreover, for $\mu\in M(f)$ we define the Lyapunov exponent of $\mu$ by 
\[
\chi(\mu):=\int \log |f'| d\mu.
\]
We set 
\[
H:=\sup\left\{\frac{h(\mu)}{\chi(\mu)}:\mu\in M(f),\chi(\mu)>0\right\}.
\]
Note that if $f$ is topologically exact and has a H\"older
continuous derivative and only non-flat critical points,
then by Ruelle's inequality (see \eqref{eq ruelle's ine})
we have
\begin{align}\label{eq upper bound H}
    H\leq 1.
\end{align}

Let $\phi$ be a continuous function on $I$. 
For all $\alpha\in \mathbb{R}$ we define 
\[
B(\alpha):=\left\{x\in I:\lim_{n\to \infty}\frac{1}{n}\sum_{k=0}^{n-1}\phi(f^k(x))=\alpha\right\}.
\]
We define the Birkhoff spectrum $b:\mathbb{R}\rightarrow [0,1]$ by 
\[
b(\alpha):=\dim_H(B(\alpha)).
\]
For each $\epsilon>0$ we define functions $D_\epsilon, D:\mathbb{R}\rightarrow \mathbb{R}\cup \{-\infty\}$ by 
\begin{align*}
&D_\epsilon(\alpha):=\sup\left\{\frac{h(\mu)}{\chi(\mu)}: \mu\in M(f), \chi(\mu)>0, \left|\int \phi d\mu-\alpha\right|<\epsilon\right\},
\\&D(\alpha):=\sup\left\{\frac{h(\mu)}{\chi(\mu)}: \mu\in M(f), \chi(\mu)>0, \int \phi d\mu=\alpha\right\}.
\end{align*}
Moreover, we define the function $D_0:\mathbb{R}\rightarrow \mathbb{R}\cup\{-\infty\}$ by 
\[
D_0(\alpha)=\lim_{\epsilon \to 0} D_\epsilon(\alpha).
\]
We set 
\[
\alpha_{\min}:=\min\left\{\int \phi d\mu: \mu\in M(f)\right\}
\text{ and } \alpha_{\max}:=\max \left\{\int \phi d\mu: \mu\in M(f)\right\}.
\]
Since $f$ is continuous, $M(f)$ is compact in the weak topology (see \cite[Theorem 6.10]{walters2000introduction}).
Consequently, the maximum and minimum in the above definition
are attained.  
Define 
\[
\hyperbolic:=\{\alpha\in [\alpha_{\min},\alpha_{\max}]: D_0(\alpha)<H\}
\text{ and }
\bad:=\{\alpha\in [\alpha_{\min},\alpha_{\max}]: D_0(\alpha)=H\}.
\]
By definition, $\bad\neq\emptyset$.
Moreover, 
$
D_0(\alpha)=
\inf_{\epsilon>0}
\sup_{|\beta-\alpha|<\epsilon}D(\beta)
$
for each $\alpha\in\mathbb {R}$.
Thus, $D_0$ is upper semicontinuous.
In particular, $\bad$ is closed.
Moreover, by Lemma \ref{lemma bad is an interval},
$\bad$ is an interval.

We are now in a position to state our main theorem.
\begin{thm}\label{thm main}
    Let $f:I\rightarrow I$ have a H\"older continuous derivative and only non-flat critical points and let $\phi$ be a continuous function. We assume that $f$ is topologically exact. Then for each $\alpha\in [\alpha_{\min},\alpha_{\max}]$ we have 
    \begin{align}\label{eq main lower bound}
    b(\alpha)\geq D_0(\alpha).
    \end{align}
    Moreover, if $\alpha\in \hyperbolic$ then we have 
    \begin{align}\label{eq conditional variational principle}
    b(\alpha)=D_0(\alpha).
    \end{align}
    Furthermore, if $\alpha\in \hyperbolic\cap (\alpha_{\min},\alpha_{\max})$ then we have 
    \begin{align}\label{eq conditional variational principle strong}
    b(\alpha)=D(\alpha).
    \end{align}
    In particular, if $H=1$ then for each $\alpha\in [\alpha_{\min},\alpha_{\max}]$ we have  
    \begin{align}\label{eq main H=1 weak}
    b(\alpha)=D_0(\alpha)
    \end{align}
    and for each $\alpha\in (\alpha_{\min},\alpha_{\max})\setminus \partial \bad$ we have 
    \begin{align}\label{eq main H=1 strong}
    b(\alpha)=D(\alpha).
    \end{align}
\end{thm}

\begin{cor}\label{cor main}
Let $f:I\to I$ be a map with a H\"older continuous
derivative and only non-flat critical points,
and let $\phi:I\to\mathbb{R}$ be continuous.
Assume that $f$ is topologically exact.
Then the function $b$ is nondecreasing on
$[\alpha_{\min},\min\bad]$ and nonincreasing on
$[\max\bad,\alpha_{\max}]$.
Moreover, $b$ is continuous on $\hyperbolic\cap(\alpha_{\min},\alpha_{\max})$.
    \end{cor}
Since $\bad$ is an interval, Corollary \ref{cor main}
implies that if $H=1$, then all points of discontinuity
of $b$ in $(\alpha_{\min},\alpha_{\max})$ lie on
the boundary of $\bad$.
In particular, if $H=1$, then $b$ has at most two
points of discontinuity in $(\alpha_{\min},\alpha_{\max})$.

For a subset $A\subset I$ satisfying $f(A)\subset A$,
we say that $f|_A$ is uniformly expanding if there
exist constants $C>0$ and $s>1$ and an integer
$N\in\mathbb{N}$ such that $|(f^n)'(x)|\geq Cs^n$
for every $x\in A$ and every $n\geq N$.

\begin{thm}\label{thm uniform compact appro spectrum}
    Let $f:I\to I$ be a map with a H\"older continuous
derivative and only non-flat critical points,
and let $\phi:I\to\mathbb{R}$ be continuous.
Assume that $f$ is topologically exact. Then there exists a sequence $\{Y_N\}_{N\in\mathbb{N}}$ of nonempty compact subsets of $I$ such that for each $N\in\mathbb{N}$ we have $Y_N\subset Y_{N+1}$, $f(Y_N)=Y_N$, $f$ is uniformly expanding on $Y_N$ and for each nonempty compact set $K\subset \hyperbolic\cap (\alpha_{\min},\alpha_{\max})$ we have 
\begin{align}\label{eq main uniform weak}
\lim_{N\to\infty}\sup_{\alpha\in K}|b(\alpha)-\dim_H(B(\alpha)\cap Y_N)|=0.
\end{align}
Moreover, if $H=1$ then for each nonempty compact set $K\subset (\alpha_{\min},\alpha_{\max})\setminus\partial\bad$ we have 
\begin{align}\label{eq main uniform strong}
\lim_{N\to\infty}\sup_{\alpha\in K}|b(\alpha)-\dim_H(B(\alpha)\cap Y_N)|=0.
\end{align}
\end{thm}

For further properties of the system $(f|_{Y_N},Y_N)$,
see Remark \ref{rem Y_N property}.

We explain how Theorems \ref{thm main} and
\ref{thm uniform compact appro spectrum},
together with Corollary \ref{cor main},
apply to a large class of $S$-unimodal maps. A continuously differentiable map $f:[0,1]\to[0,1]$
is said to be a unimodal map if $f(0)=f(1)=0$
and $f$ has a unique critical point $c\in(0,1)$.
A unimodal map $f$ is said to be an $S$-unimodal map
if $f$ is of class $C^3$, $|f'(0)|>1$, $f$ has
negative Schwarzian derivative on $[0,1]\setminus\{c\}$
(see \cite[Chapter II-6]{demelo}), and $c$ is
a non-flat critical point.
Combining Guckenheimer's classification of topological
attractors \cite{Guckenheimer}
(see also \cite[Chapter III.4]{demelo} and \cite{Bruin})
with \cite[Theorem 2.19 and Proposition 2.34]{Ruette},
we obtain the following conclusion.
If an $S$-unimodal map $f$ has no periodic attractors
and is at most finitely renormalizable, then there exist
$p\in\mathbb{N}$ and a compact interval $J$
with $c\in\operatorname{Int}(J)$ such that
$f^p(J)\subset J$, the restriction
$f^p|_J:J\to J$ is topologically exact,
and the intervals $J,f(J),\ldots,f^{p-1}(J)$
have pairwise disjoint interiors
(see also \cite{ChungRiveraTakahasi, TakahasiThujii}).
Hence, an affine rescaling of $f^p|_J$ satisfies
the assumptions of Theorem \ref{thm main} and
Corollary \ref{cor main}.
In particular, consider the family of quadratic maps
$f_a:[0,1]\to[0,1]$, $a\in(0,4]$, defined by
$f_a(x):=ax(1-x)$.
If $f_a$ has no periodic attractors and is at most
finitely renormalizable, then Theorem \ref{thm main}
and Corollary \ref{cor main} apply to an appropriate
rescaling of an iterate of $f_a$. 

Next, we give a sufficient condition for 
\[
H=1.
\]
Let $f$ be a map satisfying the assumptions of Theorem \ref{thm main}. 
For $t\in\mathbb{R}$ we define the topological pressure 
$
P(t):=\sup\left\{h(\mu)-t\chi(\mu): \mu\in M(f)\right\}.
$
By \cite[Section 1.2]{ChungRiveraTakahasi} (see also \cite[Proposition A.1]{Rivera-Letelierasymptotic}), for each $\mu\in M(f)$ we have 
\begin{align}\label{eq non negative lyapunov}
    \chi(\mu)\geq 0.
\end{align}
Moreover, we have Ruelle's inequality, that is, for each $\mu\in M(f)$ we have 
\begin{align}\label{eq ruelle's ine}
    h(\mu)\leq \chi(\mu).
\end{align}
Combining the upper bound in the full level-2 large deviation principle established in \cite[Main Theorem]{ChungRiveraTakahasi} with the nonnegativity of the rate function and its representation in \cite[(1)]{ChungRiveraTakahasi} (see also \cite[Section 2.1]{ChungRiveraTakahasi}), we obtain 
\[
P(1)=0.
\]
Thus, if we have 
$\chi_{\inf}:=\inf\{\chi(\mu):\mu\in M(f)\}>0 \text{ then }H=1.$
By Ruelle's inequality \eqref{eq ruelle's ine}, we obtain 
\[
H=\sup\{t\geq 0: P(t)> 0\}.
\]
Thus, by \cite[Proposition 9.2]{IommiTodd2010naturalequilibriumstate} (see also \cite{IommiTodd2013inducing}), if $f$ is an $S$-unimodal map with no wild attractor then we have 
\[
H=1.
\]
See \cite{BruinKellerNowickiStrienwild} for the definition
of a wild attractor and the existence of maps
with wild attractors. By \cite{Weixiao}, an $S$-unimodal map whose
critical point has order $l\in(1,2]$ has no wild
attractor (see also \cite{Lyubich}). Hence, if $f$ is an $S$-unimodal map whose
critical point has order $l\in(1,2]$, then $H=1$. In particular, if $f$ is a topologically exact quadratic map, then $H=1$.   

Let $f$ be a map satisfying the assumptions of Theorem \ref{thm main}. The arguments used in \cite{BarreiraSaussol,
JJTOPnonuniformly, MixedJT} rely heavily on
the Markov structure and the absence of critical points. However, in general, it is difficult to find a Markov
structure for $f$. Moreover, $f$ may have critical points.  
Hence, it is not straightforward to directly extend
the arguments in \cite{BarreiraSaussol,
JJTOPnonuniformly, MixedJT}. 
We study the Birkhoff spectrum of a general continuous potential $\phi$,
which does not necessarily encode the geometric
properties of $f$. Therefore, as mentioned at the beginning of the introduction,
we cannot expect the Birkhoff spectrum $b$ to be expressed
as the Legendre transform of a thermodynamic quantity.
Therefore, the methods used in
\cite{GelfertPrzytyckiRams,
IommiTodd2011dimension} to study the Lyapunov spectrum and 
the local dimension spectrum are not directly applicable
to our setting. Moreover, we do not impose the hyperbolicity
assumptions used in \cite{Chung, ChungTakahasi}.

A key ingredient in the proofs of Theorems
\ref{thm main} and \ref{thm uniform compact appro spectrum}
is a version of Katok's lemma with a uniform base set $B$
(Lemma \ref{lemma key lemma upper bound}).
We discuss this lemma in more detail below.
Chung, Rivera-Letelier and Takahasi
\cite[Lemma 2.2]{ChungRiveraTakahasi}
established a version of Katok's lemma.
However, the base set provided by their lemma
may depend on the choice of the measure $\mu$.
For this reason, we cannot directly use their lemma
to construct an induced system that simultaneously
captures the information carried by several invariant
measures.

To overcome this problem,  we first construct a nondegenerate compact interval $B$ such that for each open interval $U\subset (0,1)$ there exist $a\in\mathbb{N}$ and $A\subset U$ such that $f^a|_A:A\rightarrow B$ is a diffeomorphism (Lemma \ref{lemma outgoing conector}). Moreover, in the proof of Lemma \ref{lemma key lemma upper bound},
we construct a nondegenerate compact interval
$V\subset \operatorname{Int}(B)$ such that
$f^b|_V$ is a diffeomorphism onto its image
for some $b\in\mathbb{N}$. Combining these constructions with  a version of Katok's lemma established in \cite{ChungRiveraTakahasi}, we establish a version of Katok's lemma with a uniform base set $B$ (Lemma \ref{lemma key lemma upper bound}). Furthermore, we show that the compact intervals
$J_1,\ldots,J_M$ defining the branches can be chosen
within any prescribed nonempty open subinterval
of $\operatorname{Int}(B)$. This property plays an important role in the
construction of a sequence of induced systems
in Section \ref{sec construction}.

To prove \eqref{eq main lower bound}, using a version of Katok's lemma with a uniform base set $B$
(Lemma \ref{lemma key lemma upper bound}), for $\alpha\in [\alpha_{\min},\alpha_{\max}]$ and a suitable sequence of $f$-invariant measures, we construct an inhomogeneous Cantor set $K\subset B(\alpha)$ with large Hausdorff dimension.   

To prove Theorem \ref{thm uniform compact appro spectrum},
in Section \ref{sec construction}, we choose a suitable
sequence of invariant measures
$\{\mu_n\}_{n\in\mathbb{N}}$ with positive entropy.
Using our version of Katok's lemma with a uniform
base set $B$ (Lemma \ref{lemma key lemma upper bound}),
we construct a sequence of uniformly expanding induced
systems $\{F_N\}_{N\in\mathbb{N}}$ with full branches
such that each $F_N$ simultaneously captures the relevant
information carried by the measures $\mu_1,\ldots,\mu_N$.
We emphasize that this construction works for any
sequence of invariant measures with positive entropy,
not just for the particular sequence chosen in the proof
of Theorem \ref{thm uniform compact appro spectrum}.

To prove \eqref{eq conditional variational principle},
for each $\alpha\in[\alpha_{\min},\alpha_{\max}]$,
we construct suitable covers of $B(\alpha)\cap J$,
where $J$ is a certain nondegenerate interval (Lemma \ref{lemma nice cover}).
The idea behind this construction comes from
\cite[Section 4]{ChungRiveraTakahasi}. Moreover, we show that for each
$\alpha\in\hyperbolic$ and each
$s\in(D_0(\alpha),H)$, there exist
$\delta>0$ and $\epsilon>0$ such that
$
h(\mu)-s\chi(\mu)\leq-\epsilon
$
for every $\mu\in M(f)$ satisfying
$
\left|\int\phi\,d\mu-\alpha\right|<\delta
$
(Lemma \ref{lemma large dimension}). By applying this result to a sequence of measures constructed using the properties of the covers
established in Lemma \ref{lemma nice cover}, we establish \eqref{eq conditional variational principle}.

The paper is organized as follows.
In Section \ref{sec Auxiliary lemma 2}, we establish
auxiliary lemmas needed for the proof of
Theorems \ref{thm main} and
\ref{thm uniform compact appro spectrum}.
In Section \ref{sec uniform Katok's lemma}, we establish
a version of Katok's lemma with a uniform base set.
In Section \ref{sec proof of main lower bound}, we prove
\eqref{eq main lower bound}.
In Section \ref{sec Auxiliary lemmas 3}, we establish
auxiliary lemmas needed for the proof of
\eqref{eq conditional variational principle}.
In Section \ref{sec construction of nice covers},
we construct suitable covers.
In Section \ref{sec proof of upper bound}, we prove
\eqref{eq conditional variational principle}.
In Section \ref{sec proof of results in Sec 1},
we prove Theorem \ref{thm main} and Corollary \ref{cor main}.
In Section \ref{sec construction}, we construct
a suitable sequence of induced systems and use it
to obtain an increasing sequence of invariant
uniformly hyperbolic sets. 
In Section \ref{sec proof of Theorem unifrom}, we prove Theorem \ref{thm uniform compact appro spectrum}.

\subsection*{Notation and terminology}
Throughout this paper, in addition to the notation and terminology
introduced in Section \ref{sec Main results}, we use the following. Let $f:I\rightarrow I$ be a differentiable map and let $\psi$ be a function on $I$.
\begin{itemize}
\item For $A\subset I$, we denote by $\operatorname{Int}(A)$, $\overline{A}$ and $\partial A$ the interior, closure, and boundary of $A$,
respectively, with respect to the Euclidean topology on $\mathbb{R}$.
\item For each $x\in I$ and $r>0$ we define $B(x,r):=\{y\in I:|x-y|\leq r\}$.
\item For $A\subset \mathbb{R}$ we denote by $|A|$ the Lebesgue measure of $A$.
    \item For an interval $J\subset I$ and $l\in\mathbb{N}$ we set 
    \[
    \dist(f^l|_{J}):=\sup_{x,y\in J}\frac{|(f^l)'(x)|}{|(f^l)'(y)|}.
    \]
    \item For intervals $J,J'\subset I$ and $l\in\mathbb{N}$,
we say that $f^l|_J:J\to J'$ is a diffeomorphism
if it is bijective and
$\inf_{x\in J}|(f^l)'(x)|>0$.
\item We define the post-critical set $P(f)$ by 
\[
P(f):=\bigcup_{j=1}^\infty f^j(\text{Crit}(f)).
\]
\item For $n\in\mathbb{N}$ and $x\in I$ we set \[
S_n\psi(x):=\sum_{k=0}^{n-1}\psi(f^k(x)).\]
\item We set 
\[
\|\psi\|_{\infty}:=\sup_{x\in I}|\psi(x)|.
\]
\item For a Borel probability measure $\mu$ on $I$ we set 
\[
\mu(\psi):=\int \psi d\mu.
\]
\item For each $A\subset I$ and each Borel measurable map $F:A\rightarrow A$ we denote by $M(F)$ the set of $F$-invariant Borel probability measures on $A$.   
\end{itemize}

\section{Preliminary}\label{sec Auxiliary lemma 2}
In this section, we establish some basic properties
that will be used throughout this paper. Throughout this section,
we fix a topologically exact dynamical system
$f:I\to I$ with a H\"older continuous derivative
and only non-flat critical points, together with
a continuous function $\phi:I\to\mathbb{R}$.

\begin{lemma}\label{lemma upper semi continuous}
    The function $\mu\in M(f)\mapsto \chi(\mu)$ is upper semicontinuous.
\end{lemma}
\begin{proof}
    For each $M\in \mathbb{N}$ we define 
    \[
    \ell_M(x):=\max\{\log |f'(x)|, -M\}.
    \]
    Then the function $\ell_M$ is continuous and bounded on $I$.
    In particular, for each $M\in\mathbb{N}$ 
    \begin{align}\label{eq ell_M continuous}
        \text{the map $\mu\in M(f)\mapsto \mu(\ell_M)$ is continuous for the weak topology.}
    \end{align}
    Moreover, for all $M\in \mathbb{N}$ and $x\in I$ we have $\ell_M(x)\geq \ell_{M+1}(x)$ and $\lim_{M\to\infty}\ell_M(x)=\log |f'(x)|$. Therefore, since $\sup_{x\in I}\log |f'(x)|\leq \log \|f'\|_\infty$, the monotone convergence theorem implies that for each $\mu\in M(f)$ we have 
    \[
    \chi(\mu)=\lim_{M\to\infty} \mu(\ell_M)=\inf_{M\in\mathbb{N}}\mu(\ell_M).
    \]
    Combining this with \eqref{eq ell_M continuous}, we obtain the desired result.
\end{proof}

We denote by $h_\text{top}(f)$ the topological entropy of $f$ (see \cite[Chapter 7]{walters2000introduction} for the definition of the topological entropy). We have the following variational principle \cite[Theorem 8.6]{walters2000introduction}: 
\begin{align}\label{eq variational principle}
    h_\text{top}(f)=\sup\{h(\mu):\mu\in M(f)\}.
\end{align}
Recall that $f$ is topologically exact. 
Thus, we have the following (see, for example, \cite[Lemma 2.72]{BruinBook}): 
\begin{align}\label{eq positive topological entropy}
    h_\text{top}(f)>0.
\end{align}
Moreover, by \cite[Theorem 8.7]{Blokh} or \cite[Theorem 3.4]{Ruette}, $f$ has the specification property (see \cite[Definition 3.1]{Ruette} for the definition of the specification property). 
Therefore, we have the following properties:
\begin{itemize}
    \item[(S1)] For each $\mu\in M(f)$  the measure $\mu$ is entropy-approachable, that is, for each neighborhood $\mathcal{V}\subset M(f)$ of $\mu$ and each $h<h(\mu)$ there exists an ergodic measure $\nu\in \mathcal{V}$ such that  $h(\nu)>h$ (see, for example, \cite{EizenbergKiferWeiss} and \cite[Theorem 2.1]{PfisterSullivan}).
    \item[(S2)] For each $\mu\in M(f)$, $\mu$ has a generic point, that is, there exists $x\in I$ such that for every continuous function $\psi$ on $I$ we have $\lim_{n\to\infty}S_n\psi(x)/n=\mu(\psi)$ (see, for example, \cite[Corollary 21.15]{DenkerGrillenbergerSigmund}).
\end{itemize}
By Lemma \ref{lemma upper semi continuous}, \eqref{eq non negative lyapunov} and (S1), we obtain the following:
\begin{lemma}\label{lemma ergodic approximation}
    Let $\mu\in M(f)$. For all $\ell\in \mathbb{N}$, $\epsilon>0$ and continuous functions $\psi_1,\cdots,\psi_\ell: I\rightarrow \mathbb{R}$ there exists an ergodic measure $\nu\in M(f)$ such that $h(\nu)>h(\mu)-\epsilon$, $\chi(\nu)<\chi(\mu)+\epsilon$ and $|\mu(\psi_j)-\nu(\psi_j)|<\epsilon$ for each $1\leq j\leq \ell$.   
\end{lemma}

Furthermore, we have the following:
\begin{lemma}\label{lemma nonempty B}
    For each $\alpha\in (\alpha_{\min},\alpha_{\max})$ there exists $\mu\in M(f)$ such that $\mu(\phi)=\alpha$ and $\chi(\mu)\geq h(\mu)>0$. Moreover, for each $\alpha\in [\alpha_{\min},\alpha_{\max}]$ we have $\overline{B(\alpha)}=I$.
\end{lemma}
\begin{proof}
    Let $\alpha\in (\alpha_{\min},\alpha_{\max})$. Choose measures $\mu_1,\mu_2\in M(f)$ with $\mu_1(\phi)<\alpha<\mu_2(\phi)$. By \eqref{eq positive topological entropy} and \eqref{eq variational principle}, there exists $\nu\in M(f)$ such that $h(\nu)>0$. There exist $t\in [0,1)$ and $i\in\{1,2\}$ such that $\alpha=t\mu_i(\phi)+(1-t)\nu(\phi)$. We set $m=t\mu_i+(1-t)\nu$. Then we have $m(\phi)=\alpha$ and $m\in M(f)$. Moreover, by the affinity of the entropy map $\hat\mu\mapsto h(\hat \mu)$ (see \cite[Theorem 8.1]{walters2000introduction}) and Ruelle's inequality \eqref{eq ruelle's ine}, we obtain $\chi(m)\geq h(m)\geq (1-t)h(\nu)>0$.
    Next, we will show that for each $\alpha\in [\alpha_{\min},\alpha_{\max}]$ we have $\overline{B(\alpha)}=I$. Since $f$ is topologically exact, it is enough to show that $B(\alpha)\neq \emptyset$. Since $\alpha\in [\alpha_{\min},\alpha_{\max}]$, there exists $\nu\in M(f)$ such that $\nu(\phi)=\alpha$. Hence, the claim follows from (S2).
\end{proof}

\section{A version of Katok's lemma with a uniform base set}\label{sec uniform Katok's lemma}
In this section, we establish a version of Katok's lemma with a uniform base set. 
Throughout this section,
we fix a topologically exact dynamical system
$f:I\to I$ with a H\"older continuous derivative
and only non-flat critical points, together with
a continuous function $\phi:I\to\mathbb{R}$.

The following version of Katok's lemma was obtained in \cite[Lemma 2.2]{ChungRiveraTakahasi}:
\begin{lemma}\label{lemma Katok's lemma}
    Let $\mu\in M(f)$ be ergodic such that $h(\mu)>0$. 
    Let $l\in\mathbb{N}$ and let $\psi_1,\cdots,\psi_l:I\rightarrow \mathbb{R}$ be continuous. Then for every $\epsilon>0$ there exist integers $k\geq 2$ and $m\geq 1$, a compact interval $K\subset I$ and pairwise disjoint compact subintervals $K_1,\cdots, K_k$ of $K$ such that the following conditions hold.
    \begin{itemize}
    \item[(K1)] We have 
    \[
    \frac{\log k}{m}\geq h(\mu)-\epsilon.
    \]
        \item[(K2)] For each $1\leq i\leq k$ the map $f^m|_{K_i}:K_i\rightarrow K$ is a diffeomorphism.
        \item[(K3)] For every $x\in \bigcup_{i=1}^kK_i$ and every $j\in \{1,\cdots, l\}$ we have 
        \[
        \left|\frac{S_m\psi_j(x)}{m}-\mu(\psi_j)\right|<\epsilon \text{ and }e^{(\chi(\mu)-\epsilon)m}\leq |(f^m)'(x)|\leq e^{(\chi(\mu)+\epsilon)m}.
        \]
    \end{itemize}
\end{lemma}
\begin{rem}\label{rem coding katok}
Let $\mu\in M(f)$ be ergodic such that $h(\mu)>0$. Let $l\in\mathbb{N}$ and let $\psi_1,\cdots,\psi_l:I\rightarrow \mathbb{R}$ be continuous. Then by \eqref{eq ruelle's ine}, we have $\chi(\mu)>0$. We fix $\epsilon>0$ satisfying $\chi(\mu)>\epsilon$. Let $k\geq 2$ and $m\geq 1$ be integers, and let
$K\subset I$ be a compact interval and
$K_1,\ldots,K_k$ be pairwise disjoint compact subintervals
of $K$, as obtained in Lemma \ref{lemma Katok's lemma}. In particular, the conditions (K1), (K2) and (K3) in Lemma \ref{lemma Katok's lemma} hold.
For each $1\leq i\leq k $ we define 
\[
g_i:=(f^{m}|_{K_i})^{-1}.
\]
We set 
\[
\Sigma:=\{1,\cdots,k\}^\mathbb{N} \text{ and } \Sigma^n:=\{1,\cdots,k\}^n \text{ for each }n\in\mathbb{N}.
\]
We endow $\Sigma$
with the shift metric $d$ defined by 
$
d(\tau,\tau')=\exp(-\inf\{k\ge 1:\,\, \tau_k \neq \tau_k' \})
$ for $\tau,\tau'\in \Sigma$.
We denote by $\sigma$ the left shift map on $\Sigma$.
Since  \[
\inf _{x\in \bigcup_{i=1}^k K_i}|(f^m)'(x)|\geq e^{(\chi(\mu)-\epsilon)m}>1
\]
by (K3), the mean value theorem implies that 
\begin{align}\label{eq uniform shrinking g images}
\lim_{n\to\infty} \max_{\tau\in \Sigma^n}|g_{\tau_1}\circ\cdots\circ g_{\tau_n}(K)|=0.
\end{align}
Hence, the coding map $\pi:\Sigma\rightarrow \pi(\Sigma)$ defined by 
\[
\{\pi(\tau)\}=\bigcap_{n\in\mathbb{N}}g_{\tau_1}\circ\cdots\circ g_{\tau_n}(K)
\]
is well-defined. Moreover, since $K_1,\ldots,K_k$ are pairwise disjoint compact subintervals of $K$, the coding map is a homeomorphism and satisfies $f^m\circ \pi=\pi\circ \sigma$. Hence, the system $(f^m, \pi(\Sigma))$ is topologically conjugate to $(\sigma,\Sigma)$.
\end{rem}

\begin{lemma}
    There exist $p\in (0,1)\setminus (P(f)\cup \text{Crit}(f))$ and $l\in\mathbb{N}$ such that $f^l(p)=p$ and $|(f^l)'(p)|>1$.
\end{lemma}
\begin{proof}
By \eqref{eq positive topological entropy} and \eqref{eq variational principle}, there exists $\nu\in M(f)$ such that $h(\nu)>0$. Moreover, by (S1), there exists an ergodic measure $\mu\in M(f)$ such that $h(\mu)>0$. We fix $0<\epsilon<\chi(\mu)$. We apply Lemma \ref{lemma Katok's lemma} with the measure $\mu$ and $\epsilon$. Then by Remark \ref{rem coding katok}, there exist infinitely many points $p\in I$ such that $f^l(p)=p$ and $|(f^l)'(p)|>1$ for some $l\in\mathbb{N}$. On the other hand, since $\#\text{Crit}(f)<\infty$, the number of periodic points in $P(f)\cup \text{Crit}(f)$ is finite. This completes the proof.    
\end{proof}
We set 
\[
HP:=\{p\in (0,1)\setminus (P(f)\cup \text{Crit}(f)):\text{$f^l(p)=p$ and $|(f^l)'(p)|>1$ for some }l\in\mathbb{N}\}.
\]

\begin{lemma}\label{lemma outgoing conector}
Let $p\in HP$. There exists a compact interval $B\subset (0,1)$ such that $p\in \text{Int}(B)$ and for each open interval $U\subset (0,1)$ there exist $n\in\mathbb{N}$ and $A\subset U$  such that $f^n|_A:A\rightarrow B$ is a diffeomorphism.
\end{lemma}

\begin{proof}
    Let $p\in HP$. Then there exists $l\in\mathbb{N}$ such that $f^l(p)=p$ and $|(f^l)'(p)|>1$. Hence, there exist a compact interval $W\subset (0,1)$ and $\theta>1$ such that 
    \begin{align*}
    p\in \text{Int}(f^l(W)),\ W\subset f^l(W),\ \inf_{x\in W}|(f^l)'(x)|>\theta    
    \end{align*} 
    and $f^l|_W:W\rightarrow f^l(W)$ is a diffeomorphism. 
    We set 
    \[
    g:=(f^l|_W)^{-1} \text{ and }B:=f^l(W).
\]
For all $k\in\mathbb{N}$ we have 
\begin{align}\label{eq uniformly shrinling g^k(B)}
    |g^k(B)|\leq \theta^{-k} \text{ and }p\in g^k(B).
\end{align}
Let $U\subset (0,1)$ be an  open interval. Since $f$ is topologically exact, there exists $n\in\mathbb{N}$ such that $f^n(U)=I$. We fix $y\in U$ satisfying $f^n(y)=p$. Since $p\notin P(f)\cup \text{Crit}(f)$, we obtain $(f^n)'(y)\neq 0$. By the inverse function theorem, there exists a small open neighborhood $V\subset(0,1)$ of $y$ such that $\overline{V}\subset U$ and 
\[
f^n|_V:V\rightarrow f^n(V) \text{ is a diffeomorphism.}
\]
By \eqref{eq uniformly shrinling g^k(B)}, there exists $k\in\mathbb{N}$ such that $g^k(B)\subset f^n(V)$. Then we obtain 
\[
A:=(f^n|_V)^{-1}(g^k(B))\subset V\subset U.
\]
Moreover, $f^{kl+n}|_{A}:A\rightarrow B$ is a diffeomorphism.
\end{proof}

\begin{lemma}\label{lemma key lemma upper bound}
Let $p\in HP$ and let
$B$ be the compact interval obtained in
Lemma \ref{lemma outgoing conector}.
We fix a nonempty open interval $W\subset \text{Int}(B)$.
Let $\mu\in M(f)$ be an ergodic measure with $h(\mu)>0$,
and fix $0<\epsilon<h(\mu)$. 
Let $l\in\mathbb{N}$ and let $\psi_1,\cdots,\psi_l:I\rightarrow \mathbb{R}$ be continuous.
Then there exist $r\in\mathbb{N}$, an integer $M\geq2$ and pairwise disjoint compact subintervals $J_1,\cdots, J_M$ of $W$ such that for each $1\leq i\leq M$ the map $f^r|_{J_i}:J_i\rightarrow B$ is a diffeomorphism and for each $x\in \bigcup_{i=1}^MJ_i$ we have 
\begin{align*}
    &\frac{\log M}{r}>h(\mu)-\epsilon,\ \max_{1\leq j\leq l}\left|\frac{S_r\psi_j(x)}{r}-\mu(\psi_j)\right|<\epsilon,\ 
    \\&\left|\frac{\log |(f^r)'(x)|}{r}-\chi(\mu)\right|<\epsilon \text{ and } |(f^r)'(x)|\geq 2.
\end{align*}
\end{lemma}

\begin{proof}
    Let $0<\eta<\epsilon/4$. By applying Lemma \ref{lemma Katok's lemma} with the measure $\mu$, the small number $\eta$ and the continuous functions $\psi_1,\cdots,\psi_l$, we obtain integers $k\geq 2$ and $m\geq 1$, a compact interval $K\subset I$ and pairwise disjoint compact subintervals $K_1,\cdots, K_k$ of $K$ such that (K1), (K2) and (K3) hold.
    We use the notation introduced in Remark \ref{rem coding katok}. We fix 
    \[
    z\in \pi(\Sigma)\setminus\left(P(f)\cup \text{Crit}(f)\cup\bigcup_{n=0}^{\infty}(f^{m})^{-n}\left(\bigcup_{i=1}^k\partial K_i\right)\right).
    \]
    By applying Lemma \ref{lemma outgoing conector} with $U=\text{Int}(K)$, we obtain $a\in\mathbb{N}$ and $A\subset \text{Int}(K)$ such that 
    \begin{align}\label{eq deffio f_A^n}
    \text{$f^a|_A:A\rightarrow B$ is a diffeomorphism.} 
    \end{align}
    Since $f$ is topologically exact, there exist $b\in\mathbb{N}$ and $x\in W$ such that $f^{b}(x)=z$. Since $z\notin P(f)$, we have $(f^{b})'(x)\neq 0$. Therefore, since $z\notin \bigcup_{i=1}^k \partial K_i$, the inverse function theorem implies that there exists a compact interval $V\subset W$ such that $x\in \text{Int}(V)$, $f^b(V)\subset \text{Int}(K)$ and 
    \begin{align}\label{eq diffeo f^b_V}
    \text{$f^b|_V:V\rightarrow f^b(V)$ is a diffeomorphism.}
    \end{align}
    We set 
    \[
    K':=f^b(V).
    \]
    Then $z\in \text{Int}(K')$. Since the coding map $\pi$ is a homeomorphism, there exists $\tau\in \Sigma$ such that $z=\pi(\tau)$. Since $z\notin \bigcup_{n=0}^{\infty}(f^{m})^{-n}(\bigcup_{i=1}^k\partial K_i)$, there exists $\ell\in\mathbb{N}$ such that 
    \[
    g_{\tau|_\ell}(K)\subset \text{Int}(K'), \text{ where }g_{\tau|_\ell}:=g_{\tau_1}\circ \cdots\circ g_{\tau_\ell}.
    \]
    Recall that by \eqref{eq ruelle's ine}, we have $\epsilon<h(\mu)\leq \chi(\mu)$.
    We set 
    \[
    T:=b+\ell m+a.
    \]
    By \eqref{eq deffio f_A^n} and \eqref{eq diffeo f^b_V}, we have
    \[
    C:=
    \max\{\sup\{|\log |(f^a)'(y)||:y\in A\},\sup\{|\log |(f^b)'(y)||:y\in V\}\}<\infty.
    \]
    We choose a sufficiently large integer $N\in\mathbb{N}$
such that 
    \begin{align}\label{eq large N}
     &\frac{T(h(\mu)-\eta)}{Nm+T}<\eta,\  \frac{2C+\ell(\chi(\mu)+\epsilon)m+T\chi(\mu)}{Nm+T}<\eta,\ \frac{2T\max_{1\leq j\leq l}\|\psi_j\|_\infty}{Nm+T}<\eta 
     \nonumber\\&\text{ and } e^{(Nm+T)(\chi(\mu)-\epsilon)}\geq 2.
    \end{align}
    For each $\omega^*\in \Sigma^N$ we set $g_{\omega^*}:=g_{\omega_1^*}\circ \cdots\circ g_{\omega_N^*}$. For each $\omega^*\in \Sigma^N$ we define 
    \[
    G_{\omega^*}:=(f^b|_V)^{-1}\circ g_{\tau|_\ell}\circ g_{\omega^*}\circ (f^a|_A)^{-1}:B\rightarrow  W.
    \]
 For each $\omega^*\in \Sigma^N$ we set 
    \begin{align*}
     J_{\omega^*}:=G_{\omega^*}(B)\subset W,\  M:=k^N \text{ and }   r:=b+(N+\ell)m+a.
    \end{align*}
Relabel the intervals
$\{J_{\omega^*}\}_{\omega^*\in\Sigma^N}$
as $J_1,\ldots,J_M$.
By the definition of $G_{\omega^*}$ $(\omega^*\in \Sigma^N)$, 
    \begin{align*}
        f^{r}|_{J_{\omega^*}}:J_{\omega^*}\rightarrow B \text{ is a diffeomorphism.}
    \end{align*}
    By (K1) and \eqref{eq large N}, we have 
    \begin{align*}
        \frac{\log M}{r}\geq \frac{Nm(h(\mu)-\eta)}{Nm+T}>h(\mu)-2\eta>h(\mu)-\epsilon.
    \end{align*}
Let $x\in \bigcup_{i=1}^MJ_i$. 
Note that for each $0\leq j\leq N-1$ we have $f^{b+\ell m+jm}(x)\in \bigcup_{i=1}^k K_i$. Hence, 
by (K3) and the chain rule, we have 
\begin{align*}
    \left|\log |(f^{Nm})'(f^{b+\ell m}(x))|-Nm\chi(\mu)\right|&\leq \sum_{j=0}^{N-1}|\log |(f^{m})'(f^{b+\ell m+jm}(x))|-m\chi(\mu)|
    \\&\leq m\eta N.
\end{align*}
Since for each $0\leq j\leq \ell -1$ we have $f^{b+jm}(x)\in \bigcup_{i=1}^k K_i$, we have 
\begin{align*}
    &|\log |(f^{b+\ell m})'(x)||\leq \sup_{y\in V}|\log |(f^b)'(y)||+\sum_{j=0}^{\ell-1}|\log |(f^{m})'(f^{b+jm}(x))||
    \\&\leq\sup_{y\in V}|\log |(f^b)'(y)||+\ell(\chi(\mu)+\epsilon)m.
\end{align*}
Therefore, by \eqref{eq large N} we obtain
    \begin{align}\label{eq Lyapunov lower}
        &|\log |(f^r)'(x)|-r\chi(\mu)|\leq |\log |(f^{b+\ell m})'(x)|-(b+\ell m)\chi(\mu)|
        \\&+|\log |(f^{Nm})'(f^{b+\ell m}(x))|-Nm\chi(\mu)|
        +|\log |(f^a)'(f^{b+\ell m+Nm}(x))|-a\chi(\mu)|
        \nonumber
    \\&\nonumber\leq  2C+\ell(\chi(\mu)+\epsilon)m+T\chi(\mu) +m\eta N<2\eta r<\epsilon r.
    \end{align}
By a similar argument, we can show that for all $x\in \bigcup_{i=1}^MJ_i$ we have 
\[
\max_{1\leq j\leq l}|S_r\psi_j(x)-r \mu(\psi_j)|<\epsilon r.
\]
    By \eqref{eq Lyapunov lower} and \eqref{eq large N}, for all $x\in \bigcup_{i=1}^MJ_i$ we have
    \begin{align*}
        |(f^r)'(x)|>e^{r(\chi(\mu)-\epsilon)}\geq 2.
    \end{align*}
        \end{proof}

\section{Proof of \eqref{eq main lower bound} in Theorem \ref{thm main}}\label{sec proof of main lower bound}
In this section, we prove \eqref{eq main lower bound}
in Theorem \ref{thm main}. Throughout this section,
we fix a topologically exact dynamical system
$f:I\to I$ with a H\"older continuous derivative
and only non-flat critical points, together with
a continuous function $\phi:I\to\mathbb{R}$.
\begin{prop}\label{prop lower bound}
    For each $\alpha\in [\alpha_{\min},\alpha_{\max}]$ we have $b(\alpha)\geq D_0(\alpha)$.
\end{prop}
\begin{proof}
    If $D_0(\alpha)=0$ then there is nothing to prove. We assume that $D_0(\alpha)>0$ and fix $d>0$ satisfying $d<D_0(\alpha)$. Then there exists a sequence $\{\nu_k\}_{k\in \mathbb{N}}\subset M(f)$ such that for all $k\in\mathbb{N}$ we have $h(\nu_k)>0$, $|\nu_k(\phi)-\alpha|<k^{-1}$ and $h(\nu_k)/\chi(\nu_k)>D_0(\alpha)-k^{-1}$. Let $\{\tilde\epsilon_k\}_{k\in\mathbb{N}}$ be a sequence of positive numbers satisfying $\lim_{k\to\infty} \tilde\epsilon_k/\chi(\nu_k)=0$ and $\tilde \epsilon_k<h(\nu_k)$. By Lemma \ref{lemma ergodic approximation}, for each $k\in\mathbb{N}$ there exists an ergodic measure $\mu_k\in M(f)$ such that 
    $
     |\mu_k(\phi)-\nu_k(\phi)|<\tilde \epsilon_k 
    $, $h(\mu_k)>h(\nu_k)-\tilde \epsilon_k$ and $\chi(\mu_k)<\chi(\nu_k)+\tilde \epsilon_k$. Thus, we obtain 
    \begin{align}\label{eq estimate mu_k}
    \lim_{k\to \infty}\mu_k(\phi)=\alpha \text{ and } \liminf_{k\to\infty}\frac{h(\mu_k)}{\chi(\mu_k)}\geq\liminf_{k\to\infty}\frac{h(\nu_k)-\tilde \epsilon_k}{\chi(\nu_k)+\tilde \epsilon_k}\geq D_0(\alpha).
    \end{align}
     Let $\{\epsilon_k\}_{k\in\mathbb{N}}$ be a sequence of positive numbers satisfying 
     \begin{align}\label{eq def epsilon}
     \lim_{k\to \infty}\epsilon_k=0,\ \lim_{k\to\infty}\frac{\epsilon_k}{\chi(\mu_k)}=0 \text{ and }\epsilon_k<\frac{h(\mu_k)}{4} \text{ for each }k\in\mathbb{N}.
     \end{align}
     Let $p\in HP$ and let $B$ be the compact interval obtained in Lemma \ref{lemma outgoing conector}. By Lemma \ref{lemma key lemma upper bound}, for each $k\in\mathbb{N}$ there exist $r_k\in\mathbb{N}$, an integer $M_k\geq2$ and pairwise disjoint compact subintervals $J_{k,1},\cdots, J_{k,M_k}$ of $\text{Int}(B)$ such that for each $1\leq i\leq M_k$ the map $f^{r_k}|_{J_{k,i}}:J_{k,i}\rightarrow B$ is a diffeomorphism and for each $x\in \bigcup_{i=1}^{M_k}J_{k,i}$ we have 
\begin{align}\label{eq estimate upper bound}
    &\frac{\log M_k}{r_k}>h(\mu_k)-\epsilon_k,\ \left|\frac{S_{r_k}\phi(x)}{r_{k}}-\mu_k(\phi)\right|<\epsilon_k,\ 
    \\&\left|\frac{\log |(f^{r_k})'(x)|}{r_k}-\chi(\mu_k
)\right|<\epsilon_k \text{ and } |(f^{r_k})'(x)|\geq 2.\nonumber
\end{align}
For each $k\in \mathbb{N}$ and $1\leq i\leq M_k$
 we set 
\[
g_{k,i}:=(f^{r_k}|_{J_{k,i}})^{-1},
\]
\[
u_k:=(\chi(\mu_k)+\epsilon_k)r_k, \ \eta_k:=\epsilon_k+|\mu_k(\phi)-\alpha| \text{ and }v_k:=\log M_k.
\]
Note that by \eqref{eq estimate mu_k}, \eqref{eq def epsilon} and \eqref{eq estimate upper bound} we have 
\[
\liminf_{k\to\infty}\frac{v_k}{u_k}\geq \liminf_{k\to\infty}\frac{h(\mu_k)-\epsilon_k}{\chi(\mu_k)+\epsilon_k}\geq D_0(\alpha)>d \text{ and }\lim_{k\to\infty}\eta_k=0.
\]
Hence, by \eqref{eq estimate upper bound}, we may assume
without loss of generality that, for every $k\in\mathbb{N}$
and every $x\in\bigcup_{i=1}^{M_k}J_{k,i}$,
\begin{align}\label{eq estimate upper bound complete}
    \log |(f^{r_k})'(x)| < u_k,\ 
    \left|\frac{S_{r_k}\phi(x)}{r_k}-\alpha\right|< \eta_k,
    \text{ and }
    v_k\geq (d-\eta_k)u_k.
\end{align}

We will show that there exists a Cantor set $K$ such that 
\begin{align}
K\subset B(\alpha) \text{ and }\dim_H(K)\geq d.    
\end{align}
By \eqref{eq estimate upper bound} and \eqref{eq estimate upper bound complete}, we have $u_k>\log 2$ for all $k\in\mathbb{N}$. We inductively construct a strictly increasing sequence $\{t_k\}_{k\in\mathbb{N}}\subset \mathbb{N}$ satisfying the following conditions: for each $k\in\mathbb{N}$ we have 
\begin{align}\label{eq def t_k}
\sum_{j=1}^{k} t_j u_j\geq k^2 u_{k+1} \text{ and } \sum_{j=1}^k t_j r_j\geq k^2 r_{k+1}.
\end{align}
For each $k\in\mathbb{N}$ we set 
\[
E_{k}:=\{1,\cdots, M_k\} \text{ and }E_k^0=\{\emptyset\},
\]
where $\emptyset$ denotes the empty word. 
For each $k\in\mathbb N\cup\{0\}$
and $0\leq s<t_{k+1}$, we set
\[
E^{t_1,\ldots,t_k,s}
:=E_1^{t_1}\times\cdots\times
E_k^{t_k}\times E_{k+1}^{s},
\]
where the empty product is understood
to be a singleton.
For each $k\in\mathbb{N}$ we endow $E_k$ with the discrete topology. 
For $k\in\mathbb{N}$ we define the probability measure $\hat \mu_k^*$ on $E_k$ by setting 
\[
\hat \mu_k^*(\{\tau\})=(M_k)^{-1} \text{ for each }\tau\in E_k.
\]
Using the Kolmogorov extension theorem, we define the product measure 
\[
{\mu}^*:=\bigotimes_{i = 1}^{t_1}\hat{\mu}_1^*\times \bigotimes_{i = 1}^{t_2}\hat{\mu}_2^*\times \cdots \text{ on } \Omega:=E_1^{t_1}\times E_2^{t_2}\cdots.
\]
Here, we endow $\Omega$ with the product topology. 
For each $\ell\in\mathbb{N}$ we set
\[
T_\ell:=\sum_{q=1}^{\ell} t_q \text{ and }T_0:=0. 
\]
For each $j\in\mathbb{N}$ we set 
\[
k_j=\ell \text{ if } T_{\ell-1}<j\leq T_{\ell}.
\]
For each $k\in\mathbb{N}\cup\{0\}$, $0\leq s< t_{k+1}$ and $\tau\in E^{t_1,\cdots,t_k,s}$ we define
\[
C(\tau):=g_{k_1,\tau_1}\circ\cdots \circ 
g_{k_{|\tau|},\tau_{|\tau|}}(B),
\]
where $|\tau|:=T_k+s$.
Since $|g_{k,i}'(x)|\leq 2^{-1}$ for each $k\in\mathbb{N}$, $i\in E_k$ and $x\in B$ by \eqref{eq estimate upper bound}, for each $\omega\in \Omega$ the set 
\[
\bigcap_{n\in\mathbb{N}}C(\omega|_n), \text{ where }\omega|_n:=\omega_1\cdots\omega_n
\]
is a singleton. We define the coding map $\pi: \Omega\rightarrow  \pi(\Omega)$ by 
\[
\{\pi(\omega)\}=\bigcap_{n\in\mathbb{N}}C(\omega|_n).
\]
Since $\{J_{k,i}\}_{i\in E_k}$ is pairwise disjoint for each $k\in\mathbb{N}$, the coding map $\pi$ is a homeomorphism.
We define 
\[
K:=\pi(\Omega) \text{ and }\mu:=\mu^*\circ \pi^{-1}.
\]

We will show that 
\begin{align}\label{eq dimension k}
\dim_H(K)\geq d.    
\end{align}
By Frostman's lemma (see, for example, \cite[Theorem 15.6.3]{Urbanskinoninvertible}), it is enough to show that for all $x\in K$ we have 
\begin{align}\label{eq Frostman}
\liminf_{r\to 0}\frac{\log \mu(B(x,r)\cap K)}{\log r}\geq d.
\end{align}
For $T=T_\ell+s$ $(\ell\in \mathbb{N}\cup\{0\},\ 0\leq s< t_{\ell+1})$ we set 
\[
R_T:=\sum_{j=1}^{T} r_{k_j}, U_T:=\sum_{j=1}^{T} u_{k_j} \text{ and } V_T:=\sum_{j=1}^T v_{k_j},
\]
where $R_0=U_0=V_0=0$.
For every sufficiently small $r>0$, there exist unique
integers $\ell(r)\in\mathbb{N}$ and
$0\leq s(r)<t_{\ell(r)+1}$ such that
\[
|B|e^{-U_{T(r)+1}}<r\leq |B|e^{-U_{T(r)}},
\]
where $T(r):=T_{\ell(r)}+s(r)$.
Fix such an $r>0$. Let $x\in K$.
By the chain rule and \eqref{eq estimate upper bound complete}, for all $\omega\in \Omega$ and $\hat x\in C(\omega|_{T(r)})$ we obtain
\begin{align*}
    |(f^{R_{T(r)}})'(\hat x)|\leq e^{U_{T(r)}}.
\end{align*}
Thus, by the mean value theorem,  for all $\omega\in \Omega$ we obtain 
\begin{align*}
    |C(\omega|_{T(r)})|\geq e^{-U_{T(r)}}|B|\geq r.
\end{align*}
Hence, we obtain 
\begin{align}\label{eq ball estimate number}
    \#\{\tau\in E^{t_{1},\cdots,t_{\ell(r)},s(r)} :C(\tau)\cap B(x,r)\neq \emptyset\}\leq 4.
\end{align}
Moreover, by the definition of $\mu$, for all $\omega\in \Omega$ we have
\[
\mu(C(\omega|_{T(r)}))=e^{-V_{T(r)}}.
\]
Combining this with \eqref{eq ball estimate number}, we obtain 
\begin{align}\label{eq lower bound}
    \frac{\log \mu(B(x,r))}{\log r}\geq\frac{\log 4}{\log r}+\frac{-V_{T(r)}}{\log |B|-U_{T(r)+1}}.
\end{align}
By \eqref{eq estimate upper bound complete}, we have 
\begin{align*}
    &\liminf_{r\to 0}\frac{V_{T(r)}}{U_{T(r)}}\geq d-\limsup_{r\to 0}\frac{\sum_{j=1}^{T(r)}\eta_{k_j}u_{k_j}}{U_{T(r)}}=d \text{ and }
    \\&\limsup_{r\to 0}\frac{u_{k_{T(r)+1}}}{U_{T(r)}}\leq \limsup_{r\to 0}\frac{u_{\ell(r)+1}}{\ell(r)^2 u_{\ell(r)+1}}=0.
\end{align*}
Thus, we obtain $\liminf_{r\to 0}V_{T(r)}/U_{T(r)+1}\geq d$. Combining this with \eqref{eq lower bound}, we obtain 
\[
\liminf_{r\to 0}\frac{\log \mu(B(x,r))}{\log r}\geq d,
\]
which yields \eqref{eq Frostman}. This completes the proof of \eqref{eq dimension k}.

Finally, we show that 
\begin{align}\label{eq K is contained B}
    K\subset B(\alpha).
\end{align}
Let $x\in K$. By \eqref{eq estimate upper bound complete}, for all $T=T_\ell+s$ ($\ell\in \mathbb{N}\cup\{0\},\ 0\leq s<t_{\ell+1}$) we have 
\begin{align}\label{eq birkhoff estimate main}
    |S_{R_T}\phi(x)-R_T\alpha|\leq \sum_{j=1}^T|S_{r_{k_j}}\phi(f^{R_{j-1}}(x))-r_{k_j}\alpha| \leq \sum_{j=1}^T \eta_{k_j} r_{k_j}.
\end{align}
For every sufficiently large $n\in\mathbb{N}$, there exist unique integers $\ell(n)$ and $0\leq s(n)<t_{\ell(n)+1}$ such that 
\[
R_{T(n)}\leq n<R_{T(n)+1},
\]
where $T(n):=T_{\ell(n)}+s(n)$.
By \eqref{eq birkhoff estimate main}, we obtain 
\begin{align}\label{eq BIrkhoff finial}
    |S_n\phi(x)-n\alpha|\leq \sum_{j=1}^{T(n)}\eta_{k_j}r_{k_j}+(\|\phi\|_\infty+|\alpha|)r_{k_{T(n)+1}}.
\end{align}
By \eqref{eq def t_k}, the limit $\lim_{k\to\infty}\eta_k=0$ implies that 
\begin{align*}
\lim_{n\to\infty}\frac{\sum_{j=1}^{T(n)}\eta_{k_j}r_{k_j}}{R_{T(n)}}=0 \text{ and }
\limsup_{n\to\infty}\frac{r_{k_{T(n)+1}}}{R_{T(n)}}\leq \limsup_{n\to\infty}\frac{r_{k_{T(n)+1}}}{\ell(n)^2r_{k_{T(n)+1}}}=0.
\end{align*}
Hence, by \eqref{eq BIrkhoff finial}, we obtain $x\in B(\alpha)$. This completes the proof of \eqref{eq K is contained B}. Combining \eqref{eq K is contained B} with \eqref{eq dimension k}, we obtain the desired result.
\end{proof}

\section{Proof of \eqref{eq conditional variational principle} in Theorem \ref{thm main}}
In this section, we prove  \eqref{eq conditional variational principle} in Theorem \ref{thm main}. Throughout this section,
we fix a topologically exact dynamical system
$f:I\to I$ with a H\"older continuous derivative
and only non-flat critical points, together with
a continuous function $\phi:I\to\mathbb{R}$.

\subsection{Auxiliary lemmas}\label{sec Auxiliary lemmas 3}
The following lemma was stated in \cite[Section 3.1, p. 869]{ChungRiveraTakahasi}. However, for the convenience of the reader, we provide a proof.
\begin{lemma}\label{lemma uniform exact}
    For every $0<\kappa<1$ there exists $N(\kappa)\in\mathbb{N}$ such that for every compact interval $J\subset I$ with $|J|\geq\kappa$ we have $f^{N(\kappa)}(J)=I$.
\end{lemma}
\begin{proof}
    We define the compact set $J_\kappa$ by 
    \[
    J_\kappa:=\{(a,b)\in\mathbb{R}^2: 0\leq a\leq b\leq 1,b-a\geq \kappa\}.
    \]
    For each $(a,b)\in J_\kappa$ we fix $(u_{a,b},v_{a,b})$ with $a<u_{a,b}<v_{a,b}<b$ and set $\delta_{a,b}:=\min\{u_{a,b}-a,b-v_{a,b}\}/2$.
    We define the open neighborhood $O_{a,b}$ of $(a,b)$ by 
    \[
    O_{a,b}:=\{(a',b')\in \mathbb{R}^2: \max\{|a-a'|,|b-b'|\}<\delta_{a,b}\}.
    \]
    Since $f$ is topologically exact, for each $(a,b)\in J_\kappa$ there exists $N_{a,b}\in \mathbb{N}$ such that $f^{N_{a,b}}((u_{a,b},v_{a,b}))=I$. Since $J_\kappa$ is compact and $J_\kappa\subset \bigcup_{(a,b)\in J_\kappa} O_{a,b}$, there exist $s\in\mathbb{N}$ and $\{(a_i,b_i)\}_{1\leq i\leq s}\subset J_\kappa$ such that $J_\kappa\subset \bigcup_{i=1}^sO_{a_i,b_i}$. We set $N(\kappa):=\max_{1\leq i\leq s} N_{a_i,b_i}$. Then for each $J:=[a,b]\subset I$ with $|J|\geq \kappa$ there exists $1\leq i\leq s$ such that $(a,b)\in O_{a_i,b_i}$. Thus, we obtain $I=f^{N(\kappa)}((u_{a_i,b_i},v_{a_i,b_i})) \subset f^{N(\kappa)}(J)\subset I$.
\end{proof}

In the following, for $0<\kappa<1$ we denote by $N(\kappa)$ the natural number obtained in Lemma \ref{lemma uniform exact}.
The following lemma is an immediate consequence of Lemma \ref{lemma uniform exact}.
\begin{lemma}\label{lemma uiform shrinking}
    Let $0<\eta<1/2$. If an interval $W$ satisfies $|f^n(W)|\leq 2\eta$ then for every $0<\kappa<1$ we have $|W|<\kappa$ whenever $n\geq N(\kappa)$. 
\end{lemma}

The following lemma follows from Lemma \ref{lemma uniform exact}:
\begin{lemma}\label{lemma Birkhoff sum}
    Let $\eta\in (0,1/2)$.
    For each $n\in\mathbb{N}$ let $\mathcal{C}_n$ be the family of compact intervals $W$ satisfying $|f^n(W)|\leq 2\eta$. Then for each continuous function $\psi:I\rightarrow \mathbb{R}$ we have 
    \[
    \lim_{n\to\infty}\frac{\sup_{W\in \mathcal{C}_n}\sup_{x,y\in W}|S_n\psi(x)-S_n\psi(y)|}{n}=0.
    \]
\end{lemma}
\begin{proof}
    Let $0<\kappa<1$. Then for any $0\leq j\leq n-N(\kappa)$ and $W\in \mathcal{C}_n$ we have $|f^j(W)|<\kappa$. Therefore,  we have
    \begin{align*}
    &\limsup_{n\to\infty}\frac{\sup_{W\in \mathcal{C}_n}\sup_{x,y\in W}|S_n\psi(x)-S_n\psi(y)|}{n}
    \\&\leq \sup\{|\psi(x)-\psi(y)|:|x-y|\leq\kappa\}+\limsup_{n\to\infty}\frac{2N(\kappa)\|\psi\|_\infty}{n}
    \\&=\sup\{|\psi(x)-\psi(y)|:|x-y|\leq\kappa\}.    
    \end{align*}
    Since $\psi$ is uniformly continuous on $I$, we obtain the desired result.
\end{proof}

The following lemma will be used in the next section.
\begin{lemma}\label{lemma usefull}
    Let $0<\kappa<1$, $N\in \mathbb{N}$ and let $x_0\in (0,1)\setminus (P(f)\cup\text{Crit}(f))$. There exists $\rho^*>0$ such that for every compact interval $A\subset I$ satisfying $|A|\geq\kappa$ and $f^N(A)=I$ and every $\rho\in (0,\rho^*]$ there exists a compact interval $V\subset \text{Int}(A)$ such that $f^N|_V:V\rightarrow [x_0-2\rho,x_0+2\rho]$ is a diffeomorphism and $\dist(f^N|_{V})\leq 2$.
\end{lemma}

\begin{proof}
    Since $f^N(A)=I$ and $x_0\in (0,1)\setminus (P(f)\cup\text{Crit}(f))$, there exists $y\in \text{Int}(A)$ such that $f^N(y)=x_0$ and $(f^N)'(y)\neq 0$. By the inverse function theorem, there exists an open neighborhood $U_A$ of $y$ such that $U_A\subset \text{Int}(A)$ and $f^N|_{U_A}:U_A\rightarrow f^N(U_A)$ is a diffeomorphism. Since $x_0\in f^N(U_A)$ there exists $\epsilon_A>0$ such that $B_A:=[x_0-\epsilon_A,x_0+\epsilon_A]\subset f^N(U_A)$. 
    We set 
    \[
    V_A:=(f^N|_{U_A})^{-1}(B_A)\text{ and }d_A:=\min\{\text{dist}(V_A, \partial A),\epsilon_A\}>0.
    \]
    Note that the set 
    \[
    M_\kappa:=\{(a,b)\in\mathbb{R}^2: 0\leq a\leq b\leq 1,\  b-a\geq \kappa,\  f^{N}([a,b])=I\}
    \]
    is compact. Indeed, since the set $D:=\{(a,b)\in\mathbb{R}^2: 0\leq a\leq b\leq 1,\  b-a\geq \kappa\}$ is bounded and closed, it is enough to show that $\{(a,b)\in D:f^{N}([a,b])=I\}$ is closed. Define, for each $(a,b)\in D$,  $m(a,b):=\min _{x\in[a,b]} f^N(x)$ and $M(a,b):=\max _{x\in[a,b]} f^N(x)$. Then we have $\{(a,b)\in D:f^{N}([a,b])=I\}=\{(a,b)\in D:m(a,b)=0,\ M(a,b)=1\}$. Moreover, since $f^N$ is uniformly continuous on $[0,1]$, the functions $m$ and $M$ are continuous on $D$. Thus, $\{(a,b)\in D:f^{N}([a,b])=I\}$ is closed. 

    For each $(a,b)\in M_\kappa$ we define the open neighborhood $O_{a,b}$ of $(a,b)$ by
    \[
    O_{a,b}:=\{(a',b')\in\mathbb{R}^2:\max\{|a-a'|,|b-b'|\}<d_A/2\},
    \]
    where $A:=[a,b]$.
    Then since $M_\kappa\subset \bigcup_{(a,b)\in M_\kappa} O_{a,b}$ and $M_\kappa$ is compact, there exist $s\in\mathbb{N}$ and $\{(a_i,b_i)\}_{1\leq i\leq s}\subset M_\kappa$ such that $M_\kappa\subset \bigcup_{i=1}^s O_{a_i,b_i}$. We define 
    \[
    \tilde \rho^*:=\min _{1\leq i\leq s} \epsilon_{A_i}/2,
    \]
    where $A_i:=[a_i,b_i]$. 
    For each $1\leq i\leq s$ we set $h_i:=(f^N|_{U_{A_i}})^{-1}$. Since $f^N|_{U_{A_i}}$ is a diffeomorphism for each $1\leq i\leq s$ we have $C:=\max_{1\leq i\leq s}\max_{z\in B_{A_i}}|h_i'(z)|<\infty$. Let $\beta>0$ be such that $(f^{N})'$ is H\"older continuous with the exponent $\beta$. Then there exists $C_0>0$ such that for all $1\leq i\leq s$, $\rho\in (0,\tilde \rho^*]$ and $u,v\in [x_0-2\rho,x_0+2\rho]$ we have 
    \begin{align}\label{eq destorsion 2}
        |\log |h_i'(u)|-\log |h'_i(v)||\leq C_0|u-v|^{\beta}\leq 4^\beta C_0\rho^{\beta}.
    \end{align}
We take $\rho^*<\tilde \rho^*$ satisfying $\exp(4^\beta C_0(\rho^*)^{\beta})\leq 2$.
    
    Let $A:=[a,b]\subset I$ satisfy $|A|\geq \kappa$ and $f^N(A)=I$ and let $\rho\in (0,\rho^*]$. Since $(a,b)\in M_\kappa$ there exists $1\leq i\leq s$ such that $(a,b)\in O_{a_i,b_i}$. Note that $[x_0-2\rho,x_0+2\rho]\subset B_{A_i}$. We define $V:=h_i([x_0-2\rho,x_0+2\rho])$. Since $f^N|_{U_{A_i}}$ is a diffeomorphism, the map $f^N|_V:V\rightarrow [x_0-2\rho,x_0+2\rho]$ is a diffeomorphism and $V\subset \text{Int}(A)$. Moreover, by \eqref{eq destorsion 2}, we have $\dist(f^N|_V)\leq 2$.
    \end{proof}

\subsection{Construction of nice covers}\label{sec construction of nice covers}
In this section, we construct nice covers of
$B(\alpha)\cap[x_0-\rho,x_0+\rho]$
for every $\alpha\in[\alpha_{\min},\alpha_{\max}]$,
every $x_0\in(0,1)\setminus(P(f)\cup\operatorname{Crit}(f))$,
and every sufficiently small $\rho>0$. Our construction is inspired by \cite[Section 4]{ChungRiveraTakahasi}.

We begin by introducing the following uniform scale lemma,
which was established by Chung, Rivera-Letelier, and Takahasi
\cite[Uniform Scale Lemma]{ChungRiveraTakahasi}. 
\begin{lemma}\label{lemma uniform scale lemma}
    For every $\epsilon>0$ there exist constants $\eta_0>0$, $C>0$ and $\kappa>0$ such that for every $\eta\in (0,\eta_0)$ there exists $n_0\in\mathbb{N}$ such that the following property holds: 
    \begin{itemize}
        \item[(US)] For every $n\geq n_0$ and every subinterval $W\subset I$ with $\eta\leq |f^n(W)|\leq 2\eta$ there exist a compact subinterval $J\subset W$ and $m\in \mathbb{N}$ such that
        \begin{align*}
            &|J|\geq e^{-\epsilon n}|W|, \ n\leq m\leq n+C\log n,\ |f^m(J)|>\kappa,\ \dist(f^m|_J)\leq e^{\epsilon n}
            \\&
            \text{and }  \text{the map $f^m|_J:J\rightarrow f^m(J)$ is a diffeomorphism.}\nonumber
        \end{align*}
    \end{itemize}
\end{lemma}

Let $n\in\mathbb{N}$ and let $\eta\in (0,1/2)$. We set $M:=[\eta^{-1}]+1$. Then we have 
\begin{align}\label{eq M}
    \frac{1}{M}<\eta<\frac{3}{2M}.
\end{align}
For each $1\leq k\leq M-1$ we set 
\[
x_k=\frac{k}{M}.
\]
For each $1\leq k\leq M-1$ we denote by $\mathcal{W}_n(x_k,\eta)$ the collection of all connected components $W$ of $f^{-n}(B(x_k,\eta))$ such that $x_k\in f^n(W)$. Let 
\[
\mathcal{P}_n(\eta):=\bigcup_{k=1}^{M-1} \mathcal{W}_n(x_k,\eta).
\]
Then we have the following properties (see also \cite[Section 4]{ChungRiveraTakahasi}):
\begin{lemma}\label{lemma basic pro P}
Let $n\in\mathbb{N}$ and let $\eta\in (0,1/2)$.
We have the following properties: 
\begin{itemize}
    \item[$(\mathcal{P}1)$] For every $W\in \mathcal{P}_n(\eta)$ we have $\eta\leq |f^n(W)|\leq 2\eta$.
    \item [$(\mathcal{P}2)$] For each $x\in I$ we have $1\leq \#\{W\in \mathcal{P}_n(\eta):x\in W\}\leq 3$.
\end{itemize}     
\end{lemma}

\begin{proof}
    $(\mathcal{P}1)$ immediately follows from the construction of $\mathcal{W}_n(x_k,\eta)$.

    Let $x\in I$. We will show that $1\leq \#\{W\in \mathcal{P}_n(\eta):x\in W\}$. Since $f$ is topologically exact, we have $f^n(I)=I$. Thus, if $x\in f^{-n}(\{x_1,\cdots,x_{M-1}\})$ then $1\leq \#\{W\in \mathcal{P}_n(\eta):x\in W\}$. We assume that $x\notin f^{-n}(\{x_1,\cdots,x_{M-1}\})$. Since $\#\text{Crit}(f)<\infty$, the number of critical points of $f^n$ is finite. Hence, the set $f^{-n}(\{x_1,\cdots,x_{M-1}\})$ is finite. Let $J$ be the unique connected component of $I\setminus f^{-n}(\{x_1,\cdots,x_{M-1}\})$ such that $x\in J$. There exist $z\in \partial J$ and $1\leq k\leq M-1$ such that $f^n(z)=x_k$. Moreover, by \eqref{eq M}, we have 
    \[
    f^n([x,z])\subset B(x_k,M^{-1})\subset B(x_k,\eta).
    \]
    Hence, there exists a connected component $W$ of $f^{-n}(B(x_k,\eta))$ such that $[x,z]\subset W$. Moreover, $x_k=f^{n}(z)\in f^n(W)$. This implies that $W\in \mathcal{P}_n(\eta)$.

    Next we show that $\#\{W\in \mathcal{P}_n(\eta):x\in W\}\leq 3$.
    Notice that for each $1\leq k\leq M-1$ and each $W\in \mathcal{W}_n(x_k,\eta)$ with $x\in W$ we have $x_k\in [f^n(x)-\eta,f^n(x)+\eta]$. Moreover, since  for each $1\leq k\leq M-1$ the family $\mathcal{W}_n(x_k,\eta)$ is pairwise disjoint, we have $\#\{W\in \mathcal{W}_n(x_k,\eta): x\in W\}\leq 1$. Therefore, we obtain
    \begin{align}\label{eq count W}
        &\#\{W\in \mathcal{P}_n(\eta):x\in W\}\leq\sum_{k=1}^{M-1}\#\{W\in \mathcal{W}_n(x_k,\eta): x\in W\}
        \\&\nonumber\leq \#\{k\in\{1,\cdots,M-1\}:x_k\in  [f^n(x)-\eta,f^n(x)+\eta]\}.
    \end{align}
    By \eqref{eq M}, we have $\#\{k\in\{1,\cdots,M-1\}:x_k\in  [f^n(x)-\eta,f^n(x)+\eta]\}\leq 3$. Indeed, if there exist $4\leq j\leq M-1$ and $x_{k_1}<\cdots<x_{k_j}$ such that $x_{k_i}\in  [f^n(x)-\eta,f^n(x)+\eta]$ for each $1\leq i\leq j$ then we have 
    \[
    3M^{-1}\leq (k_{j}-k_{1})M^{-1}=x_{k_j}-x_{k_1}\leq 2\eta.
    \]
    This contradicts \eqref{eq M}. Hence, by \eqref{eq count W}, we obtain $\#\{W\in \mathcal{P}_n(\eta):x\in W\}\leq 3$.
\end{proof}

For each $\alpha\in\mathbb{R}$, $n\in\mathbb{N}$ and $\delta>0$ we set 
\[
H_n(\delta):=H_n(\alpha,\delta):=\left\{x\in I:\left|\frac{1}{n}S_n\phi(x)-\alpha\right|<\delta\right\}.
\]
For each $\eta\in(0,1/2)$ we define 
\[
    \mathcal{Q}_n(\eta,\rho,\delta):=\{W\in \mathcal{P}_n(\eta):W\cap H_n(\delta)\cap [x_0-\rho,x_0+\rho]\neq \emptyset\}.
    \]
\begin{lemma}\label{lemma nice cover}
    Let $s>0$, $\alpha\in [\alpha_{\min},\alpha_{\max}]$ and $\delta,\zeta>0$ be fixed. Let $x_0\in (0,1)\setminus (P(f)\cup\text{Crit}(f))$. There exist small numbers $\rho>0$, $0<\eta<1/2$, a constant $C>0$ and $R\in\mathbb{N}$ such that for all $n\geq R$ there exist $q\in \mathbb{N}$ with 
    \begin{align}\label{eq estimate q in the statement}
        n\leq q\leq n+C\log n
    \end{align}
       and a finite family $\{L_1,\cdots,L_N\}$ $(N\in\mathbb{N})$ of compact subintervals of $B_0:=[x_0-2\rho,x_0+2\rho]\subset (0,1)$ such that the following conditions hold:
    \begin{itemize}
        \item[(a)] We have 
        \begin{align*}
            \lim_{n\to\infty}\max\{|W|:W\in \mathcal{Q}_n(\eta,\rho,\delta)\}=0.
        \end{align*}
        \item[(b)] For each $1\leq j\leq N$ the map $f^q|_{L_j}:L_j\rightarrow B_0$ is a diffeomorphism and $\{L_1,\cdots,L_N\}$ is pairwise disjoint.
        \item[(c)] For each $1\leq j\leq N$ we have $\dist(f^q|_{L_j})\leq e^{\zeta n}$. Moreover, for each $x\in \bigcup_{j=1}^NL_j$ we have
        \begin{align*}
            \left|\frac{1}{q}S_q\phi(x)-\alpha\right|<\delta+\zeta.
        \end{align*}
        \item[(d)] We have 
        \begin{align*}
            \sum_{j=1}^N|L_j|^s\geq e^{-\zeta n}\sum_{W\in \mathcal{Q}_n(\eta,\rho,\delta)}|W|^s.
        \end{align*}
        \end{itemize}
\end{lemma}

\begin{proof}
    Let $s>0$, $\alpha\in [\alpha_{\min},\alpha_{\max}]$ and $\delta,\zeta>0$ be fixed. Let $x_0\in (0,1)\setminus (P(f)\cup\text{Crit}(f))$.
    We fix $\epsilon>0$ satisfying
    \begin{align}\label{eq epsilon}
        \max\{4s \epsilon, 2\epsilon\}< \zeta.
    \end{align}
    By Lemma \ref{lemma uniform scale lemma} and ($\mathcal{P}1$), there exist constants $\eta_0>0$, $C>0$ and $\kappa>0$ such that for every $\eta\in (0,\eta_0)$ there exists $R(\eta)\in\mathbb{N}$ such that the following property holds: 
     For every $n\geq R(\eta)$, every $\rho>0$ and every subinterval $W\in \mathcal{Q}_n(\eta,\rho,\delta)$ there exist a subinterval 
     \[
     J_W\subset W
     \]
     and $m_W\in \mathbb{N}$ such that
\begin{itemize}
    \item[(US1)] We have $|J_W|\geq e^{-\epsilon n}|W|$.
    \item[(US2)] We have $n\leq m_W\leq n+2^{-1}C\log n$.
    \item[(US3)] We have $|f^{m_W}(J_W)|>\kappa$.
    \item[(US4)] We have $\dist(f^{m_W}|_{J_W})\leq e^{\epsilon n}$.
    \item[(US5)] The map $f^{m_W}|_{J_W}:J_W\rightarrow f^{m_W}(J_W)$ is a diffeomorphism.
\end{itemize}      
    Since for all $\eta\in (0,\eta_0)$, $n\geq R(\eta)$, $\rho>0$ and $W\in \mathcal{Q}_n(\eta,\rho,\delta)$ we have $|f^{m_W}(J_W)|>\kappa$, Lemma \ref{lemma uniform exact} implies that for all $\eta\in (0,\eta_0)$, $n\geq R(\eta)$, $\rho>0$ and $W\in \mathcal{Q}_n(\eta,\rho,\delta)$ we have 
    \[
    f^{N(\kappa)+m_W}(J_W)=I.
    \]
    Moreover, by Lemma \ref{lemma usefull}, there exists $\rho(\kappa)>0$ such that for every $\eta\in (0,\eta_0)$, $n\geq R(\eta)$, $\rho\in (0,\rho(\kappa)]$ and $W\in \mathcal{Q}_n(\eta,\rho,\delta)$ there exists a compact interval \[
    V_W\subset \text{Int}(f^{m_W}(J_W))\] such that 
    \begin{align}\label{eq v_W}
    &\text{$f^{N(\kappa)}|_{V_W}:V_W\rightarrow [x_0-2\rho,x_0+2\rho]$ is a diffeomorphism 
    and}
    \\&
    \dist(f^{N(\kappa)}|_{V_W})\leq 2.    \nonumber
    \end{align}
    We fix 
    \[
    \rho\in (0,\min\{\rho(\kappa),\rho_1\}) \text{ and }\eta\in (0,\eta_0), \text{ where }\rho_1:=\min \left\{\frac{x_0}{2},\frac{1-x_0}{2}\right\}.
    \]
    We set 
    \[
{\mathcal{Q}}_n:={\mathcal{Q}}_n(\eta,\rho,\delta) \text{ for each }n\in\mathbb{N}.    
    \]
 By Lemma \ref{lemma uiform shrinking} and $(\mathcal{P}1)$, we have 
    \begin{align}\label{eq uniform shrinlinking}
        \lim_{n\to\infty} \max\{|W|: W\in\mathcal{P}_n(\eta) \}=0.
    \end{align}
    Then by \eqref{eq uniform shrinlinking}, there exists an integer $R_1\geq R(\eta)$ such that for all $n\geq R_1$ and $W\in \mathcal{Q}_n$ we have 
    \[
    W\subset B_0.
    \]
    Since $\alpha\in [\alpha_{\min},\alpha_{\max}]$, Lemma \ref{lemma nonempty B} implies that there exists $R_2\geq R_1$ such that for all $n\geq R_2$ we have 
    \begin{align}\label{eq nonempty H}
        H_n(\delta)\cap [x_0-\rho,x_0+\rho]\neq \emptyset.
    \end{align}
Set
\[
M_f:=\max\{1,\|f'\|_\infty\} \text{ and }
\Delta_n:=\frac{\sup_{W\in \mathcal{P}_n(\eta)}\sup_{x,y\in W}|S_n\phi(x)-S_n\phi (y)|}{n}.
\]
Lemma \ref{lemma Birkhoff sum} implies that there exists an integer $R\geq R_2$ such that for all $n\geq R$ we have 
\begin{align}\label{eq constant}
&\min\{|B_0|M_f^{-N(\kappa)},2^{-1}\}\geq e^{-\epsilon n} ,\  \max\left\{\Delta_n,\frac{2C\|\phi\|_\infty\log n}{n}\right\}\leq \frac{\zeta}{2},
\\&
N(\kappa)\leq 2^{-1}C\log n
\text{ and }
3([C\log n]+2)\leq e^{\epsilon sn}.\nonumber
\end{align}
Let $n\geq R$ and let $W\in {\mathcal{Q}}_n$. We define 
\begin{align}\label{eq def L_W}
L_W:=(f^{m_W}|_{J_W})^{-1}(V_W)
\end{align}
and 
\[
q_W:=m_W+N(\kappa).
\]
By (US2) and \eqref{eq constant}, we obtain
\begin{align}\label{eq estimate q_W}
n\leq q_W\leq n+2^{-1}C\log n+N(\kappa)\leq n+C\log n.
\end{align}
Note that since $f^{N(\kappa)}(V_W)=B_0$, the mean value theorem yields that  
\[
|B_0|\leq M_f^{N(\kappa)} |V_W|. 
\]
Moreover, since $f^{m_W}(L_W)=V_W$ and $L_W\subset J_W$, (US4) and the mean value theorem imply that 
\[
\frac{|V_W|}{|f^{m_W}(J_W)|}\leq e^{\epsilon n}\frac{|L_W|}{|J_W|}.
\]
Since $|f^{m_W}(J_W)|\leq 1$, (US1) implies that we have 
\begin{align}\label{eq size L_W}
|L_W|\geq e^{-\epsilon n}|J_W||V_W|\geq |B_0|M_f^{-N(\kappa)} e^{-2\epsilon n}|W|\geq e^{-3\epsilon n}|W|.
\end{align}
By the chain rule, (US4) and \eqref{eq v_W}, for all $x,y\in L_W$ we have 
\[
\frac{|(f^{q_W})'(x)|}{|(f^{q_W})'(y)|}=\frac{|(f^{N(\kappa)})'(f^{m_W}(x))||(f^{m_W})'(x)|}{|(f^{N(\kappa)})'(f^{m_W}(y))||(f^{m_W})'(y)|}\leq 2 e^{\epsilon n}\leq e^{2\epsilon n}.
\]
By \eqref{eq epsilon}, this implies that 
\begin{align}\label{eq proof tequnical lemma c dist}
    \dist(f^{q_W}|_{L_W})\leq e^{2\epsilon n}\leq e^{\zeta n}.
\end{align}
By the definition of ${\mathcal{Q}}_n$, there exists $y_W\in W\cap H_n(\delta)\cap [x_0-\rho,x_0+\rho]$. Since $L_W\subset J_W\subset W$, for all $x\in L_W$ we have 
\[
\left|\frac{S_n\phi(x)}{n}-\alpha\right|\leq \left|\frac{S_n\phi(x)}{n}-\frac{S_n\phi(y_W)}{n}\right|+\left|\frac{S_n\phi(y_W)}{n}-\alpha\right|< \Delta_n+\delta\leq \frac{\zeta}{2}+\delta.
\]
Moreover, by \eqref{eq estimate q_W}, for $x\in L_W$ we have 
\begin{align*}
&\left|\frac{S_{q_W}\phi(x)}{q_W}-\frac{S_{n}\phi(x)}{n}\right|
\leq \frac{n(q_W-n)}{q_Wn}\|\phi\|_\infty+\frac{1}{q_W}\sum_{k=n}^{q_W-1}|\phi(f^k(x))|
\\&\leq \frac{2\|\phi\|_\infty(q_W-n)}{q_W}
\leq \frac{2C\|\phi\|_\infty\log n}{n}\leq \frac{\zeta}{2}.
\end{align*}
Hence, for $x\in L_W$ we have 
\begin{align}\label{eq c birkhoff sum}
    \left|\frac{S_{q_W}\phi(x)}{q_W}-\alpha\right|< \delta+\zeta.
\end{align}

By ($\mathcal{P}2$) and \eqref{eq nonempty H}, we have ${\mathcal{Q}}_n\neq \emptyset$. 
For $n\leq q\leq n+[C\log n]+1$ we set
\[
 {\mathcal{Q}}_n(q):=\{W\in{\mathcal{Q}}_n: q_W=q \},
\]
where $[\cdot]$ denotes the floor function.
By \eqref{eq estimate q_W} we have 
\[
 {\mathcal{Q}}_n=\bigcup_{q=n}^{n+[C\log n]+1} {\mathcal{Q}}_n(q).
\]
We fix $q$ satisfying 
\begin{align}\label{eq def q}
\sum_{W\in  {\mathcal{Q}}_n(q)}|W|^s=\max_{n\leq \hat q\leq n+[C\log n]+1}\sum_{W\in  {\mathcal{Q}}_n(\hat q)}|W|^s.
\end{align}
We have 
\begin{align}\label{eq sum W}
    \sum_{W\in \mathcal{Q}_n}|W|^s=\sum_{\hat q=n}^{n+[C\log n]+1}\sum_{W\in \mathcal{Q}_n(\hat q)}|W|^s\leq ([C\log n]+2)\sum_{W\in \mathcal{Q}_n(q)}|W|^s.
\end{align}
We will show that 
\begin{align}\label{eq pairwise desjoint L_W}
  \text{Any two intervals in \(\{L_W\}_{W\in\mathcal Q_n(q)}\) are disjoint or coincide.}
\end{align}
Let $W_1,W_2\in \mathcal{Q}_n(q)$. We set 
\[
L_1:=L_{W_1},\ L_2:=L_{W_2},\ g_1:=(f^q|_{L_1})^{-1} \text{ and }g_2:=(f^q|_{L_2})^{-1}.
\]
We assume that $L_1\cap L_2\neq \emptyset$. We set 
\[
E:=\{y\in B_0:g_1(y)=g_2(y)\}.
\]
Since $g_1$ and $g_2$ are continuous on $B_0$,
the set $E$ is relatively closed in $B_0$.
Since $L_1\cap L_2\neq \emptyset$, $E$ is not empty. Moreover, the set $E$ is relatively open in $B_0$. The proof of this claim proceeds as follows: Let $y\in E$. We set $x:=g_1(y)=g_2(y)$. Since $(f^q)'(x)\neq 0$, the inverse function theorem implies that there exists an open neighborhood $U\subset I$ of $x$ such that $f^q|_U:U\rightarrow f^q(U)$ is a diffeomorphism. Since $g_1$ and $g_2$ are continuous, there exists a relatively open neighborhood $V\subset f^q(U)$ of $y$ in $B_0$ such that $g_1(V)\subset U$ and $g_2(V)\subset U$. Then for all $\hat y\in V$ we have $f^q(g_1(\hat y))=f^q(g_2(\hat y))$. Since $f^q$ is injective on $U$, we have  $g_1(\hat y)=g_2(\hat y)$. This implies that $y\in V\subset E$. Hence, $E$ is relatively open in $B_0$. Therefore, $E$ 
is a nonempty set that is both relatively closed and relatively open in $B_0$. Since $B_0$ is connected, we obtain $B_0=E$, which yields that 
\[
L_1=g_1(B_0)=g_2(B_0)=L_2.
\]
This completes the proof of \eqref{eq pairwise desjoint L_W}.

Note that if $W_1,\cdots,W_r\in \mathcal{Q}_n(q)$ satisfy $L:=L_{W_1}=\cdots=L_{W_r}$ and $W_i\neq W_j$ for $i,j\in \{1,\cdots,r\}$ with $i\neq j$ then 
\begin{align}\label{eq r upper bound}
r\leq 3.
\end{align}
Indeed, since $L_{W_j}\subset J_{W_j}\subset W_j$ for all $1\leq j\leq r$, we have $L\subset \bigcap_{j=1}^rW_{j}$. On the other hand, by ($\mathcal{P}2$) we have $\#\{W\in \mathcal{Q}_n(q):x\in W\}\leq 3$ for all $x\in L$. Thus, we obtain $r\leq 3$. We write 
\[
\{L_W\}_{W\in \mathcal{Q}_n(q)}=\{L_1,\cdots,L_N\}.
\]
Then by \eqref{eq size L_W}, \eqref{eq r upper bound}, \eqref{eq sum W} and \eqref{eq epsilon}, we obtain
\begin{align}\label{eq d}
    &\sum_{j=1}^N|L_j|^s\geq \frac{1}{3}\sum_{W\in \mathcal{Q}_n(q)}|L_W|^s\geq \frac{e^{-3\epsilon sn}}{3} \sum_{W\in \mathcal{Q}_n(q)}|W|^s\\&\geq \frac{e^{-3\epsilon sn}}{3([C\log n]+2)} \sum_{W\in \mathcal{Q}_n}|W|^s
    \geq e^{-4\epsilon sn }\sum_{W\in \mathcal{Q}_n}|W|^s\geq e^{-\zeta n}\sum_{W\in \mathcal{Q}_n}|W|^s
    . \nonumber
\end{align}
\eqref{eq estimate q in the statement} follows from \eqref{eq estimate q_W}. (a) follows from \eqref{eq uniform shrinlinking}. (b) follows from \eqref{eq v_W}, \eqref{eq def L_W} and \eqref{eq pairwise desjoint L_W}. (c) follows from \eqref{eq proof tequnical lemma c dist} and \eqref{eq c birkhoff sum}. (d) follows from \eqref{eq d}.
\end{proof}

\subsection{Upper bound of the Birkhoff spectrum}\label{sec proof of upper bound}
In this section, we give an upper bound for the Birkhoff spectrum and complete the proof of \eqref{eq conditional variational principle}.
\begin{lemma}\label{lemma large dimension}
    Let $\alpha\in \mathbb{R}$ satisfy $D_0(\alpha)<H$ and let $\max\{D_0(\alpha),0\}<s<H$. Then there exist constants $\delta,\epsilon>0$ such that for every $\mu\in M(f)$ with $|\mu(\phi)-\alpha|\leq2\delta$ we have 
    \[
    h(\mu)-s\chi(\mu)\leq -\epsilon.
    \]
    In particular, we have 
    \[
    \chi(\mu)\geq \epsilon/s>0.
    \]
\end{lemma}
\begin{proof}
    Since $s<H$, there exists $\nu\in M(f)$ such that 
    \[
C:=    h(\nu)-s\chi(\nu)>0.
    \]
    Since $D_0(\alpha)<s$, there exists $\delta>0$ such that $D_{3\delta}(\alpha)<s$. We take a sufficiently small number $\theta\in(0,1)$ so that $\theta(|\nu(\phi)-\alpha|+2\delta)<\delta$. 

Let $\mu\in M(f)$ satisfy $|\mu(\phi)-\alpha|\leq 2\delta$. We set $\tilde \mu:=(1-\theta)\mu+\theta \nu$. Then by Ruelle's inequality \eqref{eq ruelle's ine} and \eqref{eq non negative lyapunov}, we have 
\begin{align*}
    |\tilde \mu(\phi)-\alpha|\leq 2\delta+\theta(|\nu(\phi)-\alpha|+2\delta)<3\delta \text{ and }\chi(\tilde \mu)\geq \theta\chi(\nu)>0.
\end{align*}
Thus, by the affinity of the entropy map, we obtain 
\begin{align*}
(1-\theta) (h(\mu)-s\chi(\mu))+\theta C=h(\tilde \mu)-s\chi(\tilde \mu) \leq (D_{3\delta}(\alpha)-s)\chi(\tilde \mu)<0.   
\end{align*}
By setting $\epsilon:=\theta C/(1-\theta)$, we obtain the desired result.
\end{proof}
\begin{lemma}\label{lemma nice measures}
Let $q,N\in\mathbb{N}$ and let $B_0\subset (0,1)$ be a nondegenerate compact interval. Let $\alpha\in \mathbb{R}$ and let $\{L_{j}\}_{1\leq j\leq N}$ be a family of pairwise disjoint compact subintervals of $B_0$ such that for each $1\leq j\leq N$ the map $f^q|_{L_j}:L_j\rightarrow B_0$ is a diffeomorphism. We assume that there exists $\delta>0$ such that for all $x\in \bigcup_{j=1}^N L_j$ we have $|S_q\phi(x)/q-\alpha|<\delta$. Then for any $s>0$ there exists $\mu\in M(f)$ such that $|\mu(\phi)-\alpha|<\delta$ and
\[
q(h(\mu)-s\chi(\mu))\geq \log \sum_{j=1}^N|L_j|^s-s\log |B_0|-s\log \max_{1\leq j\leq N} \dist(f^q|_{L_j}).
\]
\end{lemma}

\begin{proof}
    Define $\Sigma:=\{1,\cdots,N\}^\mathbb{N}$ and $\Sigma^n:=\{1,\cdots,N\}^n$ ($n\in\mathbb{N}$). We denote by $\sigma$ the left shift map on $\Sigma$. For each $1\leq j\leq N$ we set 
    $g_j:=(f^q|_{L_j})^{-1}$. By Lemma \ref{lemma uniform exact}, we have 
    \[
    \lim_{n\to\infty} \max_{\tau\in \Sigma^n}|g_{\tau_1}\circ\cdots\circ g_{\tau_n}(B_0)|=0.
    \]
    Thus, we can define the coding map $\pi:\Sigma\rightarrow \pi(\Sigma)$ by 
    \[
    \{\pi(\tau)\}=\bigcap_{n\in\mathbb{N}}g_{\tau_1}\circ\cdots\circ g_{\tau_n}(B_0).
    \]
    Since $\{L_{j}\}_{1\leq j\leq N}$ is pairwise disjoint, the coding map $\pi$ is a homeomorphism and $f^q\circ \pi=\pi\circ \sigma$.
    For each $1\leq j\leq N$ we set 
    \[
    p_j=\frac{|L_j|^s}{\sum_{i=1}^N|L_i|^s}.
    \]
    Let $\nu^*$ be the Bernoulli measure on
$\Sigma$ associated with the probability
vector $\{p_j\}_{j=1}^N$.
    We define $\nu=\nu^*\circ \pi^{-1}$.
    Then $\nu$ is $f^q$-invariant and $h(\nu)=-\sum_{j=1}^Np_j\log p_j$, where $h(\nu)$ denotes the measure-theoretic entropy of $\nu$ with respect to the dynamical system $(f^q,\pi(\Sigma))$. By the mean value theorem, for all $1\leq j \leq N $ and $x\in L_j$ we have 
\[ 
\log|(f^q)'(x)|\leq \log |B_0|-\log |L_j|+\log \dist(f^q|_{L_j}).
\]
Therefore, we obtain 
\begin{align*}
\nonumber
   & h(\nu)-s\nu(\log|(f^q)'| )\geq -\sum_{j=1}^N p_j\log p_j-s \sum_{j=1}^Np_j \log \frac{|B_0|}{|L_j|}-s\log\max_{1\leq j\leq N} \dist(f^q|_{L_j})
   \\&=\log \sum_{j=1}^N|L_j|^s-s\log |B_0|-s\log\max_{1\leq j\leq N} \dist(f^q|_{L_j}) \text{ and }\left|\nu(\frac{S_q\phi}{q})-\alpha\right|<\delta.
\end{align*}
Therefore, the desired result follows from Abramov's
and Kac's formulas. More precisely, we have
\[
h(\nu)=qh(\mu),\ 
\nu(\log |(f^q)'|)=q\chi(\mu),
\ 
\nu(S_q\phi)=q\mu(\phi),
\]
where $\mu:=q^{-1}\sum_{t=0}^{q-1}\nu\circ f^{-t}$.
\end{proof}

\begin{prop}\label{prop upper bound}
    For each $\alpha\in [\alpha_{\min},\alpha_{\max}]$ with $D_0(\alpha)<H$ we have $b(\alpha)\leq D_0(\alpha)$.
\end{prop}
\begin{proof}
Let $D_0(\alpha)<s<H$. 
By Lemma \ref{lemma nonempty B}, $\alpha\in [\alpha_{\min},\alpha_{\max}]$ implies that $D_0(\alpha)\geq 0$.
By Lemma \ref{lemma large dimension}, there exist $\delta,\epsilon>0$ such that for every $\mu\in M(f)$ with $|\mu(\phi)-\alpha|\leq 2\delta$ we have 
\begin{align}\label{eq proof upper bound 1}
h(\mu)-s\chi(\mu)\leq -\epsilon.
\end{align}
We take a sufficiently small number $\zeta>0$ so that
\[
\zeta<\delta \text{ and }(1+s)\zeta<\epsilon/4.
\]
We fix $x_0\in (0,1)\setminus (P(f)\cup \text{Crit}(f))$. 
By Lemma \ref{lemma nice cover}, there exist small numbers $\rho>0$, $0<\eta<1/2$, a constant $C>0$ and $R\in\mathbb{N}$  such that for all $n\geq R$ there exist $q\in \mathbb{N}$ with 
    $
        n\leq q\leq n+C\log n
    $
       and a finite family $\{L_1,\cdots,L_N\}$ $(N\in\mathbb{N})$ of compact subintervals of $B_0:=[x_0-2\rho,x_0+2\rho]\subset (0,1)$ such that the conditions (a), (b), (c) and (d) in Lemma \ref{lemma nice cover} hold. Hence, by Lemma \ref{lemma nice measures}, for each $n\geq R$ there exists $\mu_n\in M(f)$ such that $|\mu_n(\phi)-\alpha|<2\delta$ and 
       \[
q(h(\mu_n)-s\chi(\mu_n))\geq \log \sum_{j=1}^N|L_j|^s-s\log |B_0|-s\log \max_{1\leq j\leq N} \dist(f^q|_{L_j}).
\]
For each $n\in\mathbb{N}$ we set
$\mathcal{Q}_n:=\mathcal{Q}_n(\eta,\rho,\delta)$.
Therefore, by \eqref{eq proof upper bound 1} and (c) and (d) in Lemma \ref{lemma nice cover}, for all $n\geq R$ we have 
\begin{align}\label{eq estimate Ws}
\sum_{W\in \mathcal{Q}_n}|W|^s\leq e^{\zeta n}\sum_{j=1}^N|L_j|^s\leq e^{(s+1)\zeta n}e^{-\epsilon q}|B_0|^s\leq e^{-\epsilon n/2}.
\end{align}
For each $\ell\in \mathbb{N}$ we set 
\[
F_\ell(\delta):=\left\{x\in [x_0-\rho,x_0+\rho]: \left|\frac{1}{k}S_k\phi(x)-\alpha\right|< \delta \text{ for each }k\geq \ell\right\}.
\]
Then by ($\mathcal{P}2$) in Lemma \ref{lemma basic pro P} for each $\ell\geq R$ and $n\geq \ell$ we have 
\[
F_\ell(\delta)\subset H_n(\delta)\cap[x_0-\rho,x_0+\rho]\subset\bigcup_{W\in \mathcal{Q}_n} W. 
\]
Hence, by (a) in Lemma \ref{lemma nice cover} and \eqref{eq estimate Ws},  for each $\ell\geq R$  we obtain $\dim_H(F_\ell(\delta))\leq s$. Since
$
B(\alpha)\cap [x_0-\rho,x_0+\rho] \subset \bigcup_{\ell=R}^\infty F_\ell(\delta)$, \cite[Theorems 15.3.4 and 15.3.5]{Urbanskinoninvertible} imply that 
\begin{align}\label{eq local}
\dim_H(B(\alpha)\cap [x_0-\rho,x_0+\rho])\leq s.    
\end{align}
Since $f$ is topologically exact, there exists $m\in\mathbb{N}$ such that $f^m((x_0-\rho,x_0+\rho))=I$. The boundedness of $\phi$ implies that $B(\alpha)\subset f^m(B(\alpha)\cap [x_0-\rho,x_0+\rho])$. Hence, by \eqref{eq local}, we obtain $b(\alpha)\leq s$. 
Since $s$ is an arbitrary number satisfying $D_0(\alpha)<s<H$, we obtain $b(\alpha)\leq D_0(\alpha)$.  
\end{proof}

Combining Proposition \ref{prop upper bound} with Proposition \ref{prop lower bound}, we obtain the following:
\begin{prop}\label{prop weak conditional}
    For each $\alpha\in [\alpha_{\min},\alpha_{\max}]$ with $D_0(\alpha)<H$ we have $b(\alpha)=D_0(\alpha)$.
\end{prop}

\section{Proofs of Theorem \ref{thm main} and Corollary \ref{cor main}}\label{sec proof of results in Sec 1}
In this section, we complete the proofs of Theorem \ref{thm main} and Corollary \ref{cor main}.  Throughout this section,
we fix a topologically exact dynamical system
$f:I\to I$ with a H\"older continuous derivative
and only non-flat critical points, together with
a continuous function $\phi:I\to\mathbb{R}$.
\begin{lemma}\label{lemma eqality}
    For $\alpha\in (\alpha_{\min},\alpha_{\max})$ with $D_0(\alpha)<H$ we have $D_0(\alpha)=D(\alpha)$.
\end{lemma}
\begin{proof}
    Let $D_0(\alpha)<s<H$.  
    We take measures $\mu_1,\mu_2\in M(f)$ satisfying $\mu_1(\phi)<\alpha<\mu_2(\phi)$.
    By Lemma \ref{lemma nonempty B} and Lemma \ref{lemma large dimension}, there exist $\epsilon,\delta>0$ such that for all $\mu\in M(f)$ with $|\mu(\phi)-\alpha|<\delta$  we have 
    \begin{align}\label{eq positive lyapunov}
        \chi(\mu)\geq c:=\epsilon/s.
    \end{align}
    By the definition of $D_0(\alpha)$, there exists a sequence $\{\nu_n\}_{n\in\mathbb{N}}\subset M(f)$ such that for all $n\in\mathbb{N}$ we have $\chi(\nu_n)>0$ and $|\nu_n(\phi)-\alpha|<\min\{\delta,|\mu_1(\phi)-\alpha|,|\mu_2(\phi)-\alpha|\}$ and 
    \[
    D_0(\alpha)=\lim_{n\to\infty} \frac{h(\nu_n)}{\chi(\nu_n)} \text{ and } \lim_{n\to\infty} \nu_n(\phi)=\alpha.
    \]
    Moreover, for all $n\in\mathbb{N}$ there exist $\mu_n\in \{\mu_1,\mu_2\}$ and $t_n\in (0,1]$ such that $\alpha=t_n\nu_n(\phi)+(1-t_n)\mu_n(\phi)$. Since $\lim_{n\to\infty}\nu_n(\phi)=\alpha$, we obtain $\lim_{n\to\infty}t_n=1$. For all $n\in\mathbb{N}$ we set $m_n:=t_n\nu_n+(1-t_n)\mu_n$. Then by \eqref{eq positive lyapunov}, for all $n\in\mathbb{N}$ we have $m_n\in M(f)$, $m_n(\phi)=\alpha$ and $\chi(m_n)\geq c>0$.  By the affinity of the entropy map, Ruelle's inequality \eqref{eq ruelle's ine} and \eqref{eq positive lyapunov}, for all $n\in\mathbb{N}$ we have 
    \[
    \left|\frac{h(m_n)}{\chi(m_n)}-\frac{h(\nu_n)}{\chi(\nu_n)}\right|=\frac{(1-t_n)|\chi(\nu_n)h(\mu_n)-\chi(\mu_n)h(\nu_n)|}{\chi(m_n)\chi(\nu_n)}\leq \frac{2(\log \|f'\|_{\infty})^2 (1-t_n)}{c^2}.
    \]
    This implies that 
    \[
    D_0(\alpha)=\lim_{n\to\infty}\frac{h(\nu_n)}{\chi(\nu_n)}=\lim_{n\to\infty}\frac{h(m_n)}{\chi(m_n)}\leq D(\alpha).
    \]
    By the definitions of $D_0(\alpha)$ and $D(\alpha)$, we have $D(\alpha)\leq D_0(\alpha)$. Thus, we are done.
\end{proof}

\begin{lemma}\label{lemma bad is an interval}
    $\bad$ is an interval.
\end{lemma}

\begin{proof}
    Let $\min \bad\leq \alpha\leq \max \bad$. We show that $\alpha\in \bad$. If $\alpha\in \{\max \bad, \min \bad\}$ then the claim is clear. Thus, we assume that $\min \bad< \alpha< \max \bad$. Let $\epsilon>0$ satisfy $H-\epsilon>0$. By Lemma \ref{lemma nonempty B}, there exist $\mu_1,\mu_2\in M(f)$ such that $\mu_1(\phi)<\alpha<\mu_2(\phi)$, $\chi(\mu_1)>0$, $\chi(\mu_2)>0$, $h(\mu_1)/\chi(\mu_1)>D_0(\min \bad)-\epsilon$ and $h(\mu_2)/\chi(\mu_2)>D_0(\max \bad)-\epsilon$. Then there exists $t\in (0,1)$ such that $\alpha=t \mu_1(\phi)+(1-t)\mu_2(\phi)$. We set $\nu:=t\mu_1+(1-t)\mu_2$. Then we have $\chi(\nu)>0$, $\alpha=\nu(\phi)$ and 
    \[
   H \geq D_0(\alpha)\geq \frac{h(\nu)}{\chi(\nu)} =\frac{th(\mu_1)+(1-t)h(\mu_2)}{t\chi(\mu_1)+(1-t)\chi(\mu_2)}\geq \min\left\{\frac{h(\mu_1)}{\chi(\mu_1)},\frac{h(\mu_2)}{\chi(\mu_2)}\right\}\geq H-\epsilon.
    \]
    Letting $\epsilon\to 0$, we obtain $\alpha\in \bad$.
\end{proof}

\emph{Proof of Theorem \ref{thm main}.} \eqref{eq main lower bound} follows from Proposition \ref{prop lower bound}. \eqref{eq conditional variational principle} follows from Proposition \ref{prop weak conditional}. \eqref{eq conditional variational principle strong} follows from Proposition \ref{prop weak conditional} and Lemma \ref{lemma eqality}.

We assume that $H=1$. \eqref{eq main H=1 weak} follows from \eqref{eq main lower bound} and \eqref{eq conditional variational principle}. To prove \eqref{eq main H=1 strong}, it suffices to show that if $\text{Int}(\bad)\neq \emptyset$ then for each $\alpha\in \text{Int}(\bad)$ we have 
\[
D(\alpha)\geq H.
\]
We assume that $\text{Int}(\bad)\neq \emptyset$. Let $\alpha\in \text{Int}(\bad)$ and let $\epsilon>0$. Since $\bad$ is a closed interval, there exist $\alpha_-,\alpha_+\in \bad$ such that $\alpha_-<\alpha<\alpha_+$. Since $D_0(\alpha_-)=D_0(\alpha_+)=H$, there exist $\mu_-,\mu_+\in M(f)$ such that $\chi(\mu_-)>0$, $\chi(\mu_+)>0$, $h(\mu_-)/\chi(\mu_-)>H-\epsilon$, $h(\mu_+)/\chi(\mu_+)>H-\epsilon$ and $\mu_-(\phi)<\alpha <\mu_+(\phi)$. Then there exists $t\in (0,1)$ such that $\alpha=t\mu_{-}(\phi)+(1-t)\mu_{+}(\phi)$. Hence, we obtain 
\[
D(\alpha)\geq \frac{h(\nu)}{\chi(\nu)}=\frac{th(\mu_-)+(1-t)h(\mu_+)}{t\chi(\mu_-)+(1-t)\chi(\mu_+)}\geq H-\epsilon,
\]
where $\nu=t\mu_-+(1-t)\mu_+$. Letting $\epsilon\to 0$, we obtain $D(\alpha)\geq H$.
\qed

\emph{Proof of Corollary \ref{cor main}.}
We first show that the function $b$ is nondecreasing on
$[\alpha_{\min},\min\bad]$. By Theorem \ref{thm main}, for all $\alpha\in [\alpha_{\min},\min\bad)$ we have $b(\alpha)=D_0(\alpha)<H\leq b(\min \bad)$. Hence, it is enough to show that for each pair $\alpha_1,\alpha_2\in [\alpha_{\min},\min\bad)$ with $\alpha_1<\alpha_2$ we have $b(\alpha_1)\leq b(\alpha_2)$. Let $\epsilon>0$. Note that $\alpha_{\min}\leq \alpha_1<\alpha_2<\min \bad$ and $H=D_0(\min \bad)\leq b(\min \bad)$ by Theorem \ref{thm main}. By Lemma \ref{lemma nonempty B} and Theorem \ref{thm main}, there exist $\mu_1,\mu_2\in M(f)$ such that $\mu_1(\phi)<\alpha_2<\mu_2(\phi)$, $\chi(\mu_1)>0$, $\chi (\mu_2)>0$, $b(\alpha_1)-\epsilon<h(\mu_1)/\chi(\mu_1)$ and $b(\alpha_1)-\epsilon< H-\epsilon<h(\mu_2)/\chi(\mu_2)$. Since $\mu_1(\phi)<\alpha_2<\mu_2(\phi)$, there exists $t\in (0,1)$ such that $\alpha_2=t\mu_1(\phi)+(1-t)\mu_2(\phi)$. We set $\nu:=t\mu_1+(1-t)\mu_2$. Then we have $\nu(\phi)=\alpha_2$. By Theorem \ref{thm main}, we obtain 
\[
b(\alpha_2)\geq \frac{h(\nu)}{\chi(\nu)}=\frac{th(\mu_1)+(1-t)h(\mu_2)}{t\chi(\mu_1)+(1-t)\chi(\mu_2)}\geq \min\left\{\frac{h(\mu_1)}{\chi(\mu_1)},\frac{h(\mu_2)}{\chi(\mu_2)}\right\}\geq b(\alpha_1)-\epsilon.
\]
Letting $\epsilon\to 0$, we obtain $b(\alpha_2)\geq b(\alpha_1)$. By a similar argument, we can show that $b$ is nonincreasing on
$[\max\bad,\alpha_{\max}]$.

Next, we show that $b$ is continuous on $\hyperbolic\cap(\alpha_{\min},\alpha_{\max})$. To this end, we will show that $D$ is lower semicontinuous on $(\alpha_{\min},\alpha_{\max})$. Let $\alpha\in (\alpha_{\min},\alpha_{\max})$. By Lemma \ref{lemma nonempty B}, there exist $\mu_1,\mu_2\in M(f)$ such that $\chi(\mu_1)>0$, $\chi(\mu_2)>0$ and $\mu_1(\phi)<\alpha<\mu_2(\phi)$. 
Let $\epsilon>0$. By Lemma \ref{lemma nonempty B}, there exists $\mu\in M(f)$ such that $\chi(\mu)>0$, $\mu(\phi)=\alpha$ and $h(\mu)/\chi(\mu)>D(\alpha)-\epsilon$.
For each $\beta\in (\mu_1(\phi), \mu_2(\phi))$ there exist $t_\beta\in [0,1]$ and $i_\beta\in \{1,2\}$ such that $\beta=t_\beta\mu(\phi)+(1-t_\beta)\mu_{i_\beta}(\phi)$. Then we have $\lim_{\beta\to \alpha} t_{\beta}=1$. For each $\beta\in (\mu_1(\phi), \mu_2(\phi))$ we set $\nu_\beta:=t_\beta\mu+(1-t_\beta)\mu_{i_\beta}$. By the definition of $D$, we have 
\begin{align*}
    \liminf _{\beta\to \alpha} D(\beta)\geq \liminf _{\beta\to \alpha}\frac{h(\nu_\beta)}{\chi(\nu_\beta)}=\liminf _{\beta\to \alpha}\frac{t_\beta h(\mu)+(1-t_\beta)h(\mu_{i_\beta})}{t_\beta \chi(\mu)+(1-t_\beta)\chi(\mu_{i_\beta})}=\frac{h(\mu)}{\chi(\mu)}>D(\alpha)-\epsilon. 
\end{align*} 
Letting $\epsilon\to 0$, we obtain $\liminf _{\beta\to \alpha} D(\beta)\geq D(\alpha)$. This implies that $D$ is lower semicontinuous on $(\alpha_{\min},\alpha_{\max})$.

By Theorem \ref{thm main}, for all $\alpha\in \hyperbolic \cap (\alpha_{\min},\alpha_{\max})$ we have $b(\alpha)=D_0(\alpha)=D(\alpha)$. Moreover, the function $D_0$ is upper semicontinuous on $\mathbb{R}$ and $D$ is lower semicontinuous on $(\alpha_{\min},\alpha_{\max})$. Hence, $b$ is continuous on $\hyperbolic \cap (\alpha_{\min},\alpha_{\max})$. 
\qed

\section{Proof of Theorem \ref{thm uniform compact appro spectrum}}\label{sec proof of Theorem unifrom}
In this section, we show Theorem \ref{thm uniform compact appro spectrum}. Throughout this section,
we fix a topologically exact dynamical system
$f:I\to I$ with a H\"older continuous derivative
and only non-flat critical points, together with
a continuous function $\phi:I\to\mathbb{R}$.
Moreover, we fix $p\in HP$ and the compact interval $B$ obtained in
Lemma \ref{lemma outgoing conector}.

\subsection{Construction of sequences of induced systems
and uniformly hyperbolic sets}\label{sec construction}
We denote by $\mathcal{UB}$ the closed unit ball of the space of continuous functions on $I$ with respect to the supremum norm. In the following, we fix a countable dense subset $\{\psi_j\}_{j\in\mathbb{N}}$ of $\mathcal{UB}$. 
We endow $M(f)$ with the weak topology. Since $I$ is a compact metric space, $M(f)$ is compact and
metrizable (see \cite[Chapter 6]{walters2000introduction}). In particular, $M(f)$ is second countable. Hence, $M(f)\times \mathbb{R}^2$ is second countable. Define
\[
\mathcal{V}:=\{(\mu,h(\mu),\chi(\mu)):
\mu\in M(f)\text{ is ergodic and }h(\mu)>0\}.
\]
Since $\mathcal{V}$ is a subspace of the second countable
space $M(f)\times\mathbb{R}^2$, there exists
$\{(\mu_i,h(\mu_i),\chi(\mu_i))\}_{i\in\mathbb{N}}\subset\mathcal{V}$ such that, for every $N\in\mathbb{N}$,
the set
$\{(\mu_i,h(\mu_i),\chi(\mu_i))\}_{i=N}^{\infty}$
is dense in $\mathcal{V}$.
We fix a family $\{W_{i}\}_{i\in\mathbb{N}}$ of nonempty open intervals such that for each $i,j\in\mathbb{N}$ with $i\neq j$ we have $W_i\cap W_j=\emptyset$ and $W_i\subset \text{Int}(B)$.
We also fix a sequence $\{\epsilon_i\}_{i\in\mathbb{N}}\subset (0,\infty)$ satisfying 
\begin{align}\label{eq epsilon i}
    \epsilon_i <\min\{2^{-i}, h(\mu_i)/2\} \text{ for each }i\in\mathbb{N}.
\end{align}
For each $i\in\mathbb{N}$, 
applying Lemma \ref{lemma key lemma upper bound} with $W_i$, $\mu_i$, $\epsilon_i$ and $\psi_1,\cdots,\psi_i$, we obtain $r_i\in\mathbb{N}$, an integer $M_i\geq 2$ and a family $\mathcal{J}_i:=\{J_{i,1},\cdots, J_{i,M_i}\}$ of pairwise disjoint compact subintervals of $W_i$ such that
 for each $1\leq j\leq M_i$ the map $f^{r_i}|_{J_{i,j}}:J_{i,j}\rightarrow B$ is a diffeomorphism and for each $x\in \bigcup_{j=1}^{M_i}J_{i,j}$ we have 
\begin{align}\label{eq apploximation cylinder i}
    &\frac{\log M_i}{r_i}>h(\mu_i)-\epsilon_i,\ \max_{1\leq k\leq i}\left|\frac{S_{r_i}\psi_{k}(x)}{r_i}-\mu_i(\psi_{k})\right|<\epsilon_i,\ 
    \\&\left|\frac{\log |(f^{r_i})'(x)|}{r_i}-\chi(\mu_i)\right|<\epsilon_i \text{ and } |(f^{r_i})'(x)|\geq 2.\nonumber
\end{align}

\begin{lemma}\label{lemma any measure aproximated for some i}
    Let $\nu\in M(f)$ satisfy $h(\nu)>0$. Then for each $\delta>0$ and $N\in\mathbb{N}$ there exists $i\geq N$ such that for all $x\in \bigcup_{j=1}^{M_i}J_{i,j}$ we have 
    \[
    \frac{\log M_i}{r_i}>h(\nu)-\delta, \ \frac{1}{r_i}\log |(f^{r_i})'(x)|<\chi(\nu)+\delta \text{ and } \left|\frac{S_{r_i}\phi(x)}{r_i}-\nu(\phi)\right|<\delta.
    \]
\end{lemma}
\begin{proof}
    We first show that 
    \begin{align}\label{eq i 1}
        \lim_{i\to\infty}\sup_{x\in \bigcup_{j=1}^{M_i}J_{i,j}}\left|\frac{1}{r_i}S_{r_i}\phi(x)-\mu_i(\phi)\right|=0.
    \end{align}
    If $\phi\equiv 0$ then there is nothing to prove. Thus, we assume that $\|\phi\|_\infty>0$. We fix $\eta>0$. Since $\{\psi_j\}_{j\in\mathbb{N}}$ is a dense subset of $\mathcal{UB}$ and $\lim_{i\to\infty}\epsilon_i=0$, there exists $k\in\mathbb{N}$ such that 
    \[
    \left\|\frac{\phi}{\|\phi\|_\infty}-\psi_k\right\|_\infty<\frac{\eta}{4\|\phi\|_\infty} \text{ and }\|\phi\|_{\infty}\epsilon_i<\frac{\eta}{2} \text{ for each }i\geq k.
    \]
    By \eqref{eq apploximation cylinder i}, for each $i\geq k$ and $x\in \bigcup_{j=1}^{M_i}J_{i,j}$ we have 
    \begin{align*}
        &\left|\frac{S_{r_i}\phi(x)}{r_i}-\mu_i(\phi)\right|
        \leq\frac{|S_{r_i}(\phi-\|\phi\|_\infty\psi_k)(x)|}{r_i} + \|\phi\|_\infty\left|\frac{S_{r_i}\psi_k(x)}{r_i}-\mu_i(\psi_k)\right|
        \\&+|\mu_i(\|\phi\|_\infty\psi_k-\phi)|<\eta/4+\|\phi\|_\infty \epsilon_i+\eta/4<\eta.
    \end{align*}
    This completes the proof of \eqref{eq i 1}.
    We fix $\tilde N\geq N$ such that for each $i\geq \tilde N$ we have 
    \begin{align}\label{eq tilde N}
        \sup_{x\in \bigcup_{j=1}^{M_i}J_{i,j}}\left|\frac{1}{r_i}S_{r_i}\phi(x)-\mu_i(\phi)\right|<\frac{\delta}{3} \text{ and }\epsilon_i<\frac{\delta}{3}.
    \end{align}

    By Lemma \ref{lemma ergodic approximation}, there exists an ergodic measure $m\in M(f)$ such that $h(m)>0$, $h(m)>h(\nu)-\delta/3$, $\chi(m)<\chi(\nu)+\delta/3$ and $|m(\phi)-\nu(\phi)|<\delta/3$. Since the map $\mu\mapsto \mu(\phi)$ is continuous, the set 
    \[
    \mathcal{U}:=\{(\mu,a,b)\in M(f)\times \mathbb{R}^2:\max\{ |\mu(\phi)-m(\phi)|,|a-h(m)|,|b-\chi(m)|\}<\delta/3\}
    \]
    is open and contains $(m,h(m),\chi(m))$. Hence,  $\mathcal{U}\cap \mathcal{V}\neq \emptyset$. Therefore,  since the set 
$\{(\mu_i,h(\mu_i),\chi(\mu_i))\}_{i=\tilde N}^{\infty}$
    is dense in $\mathcal{V}$, there exists $i\geq \tilde N\ge N$ such that 
    \[
    (\mu_i,h(\mu_i),\chi(\mu_i))\in \mathcal{U}.
    \]
    Thus, by \eqref{eq apploximation cylinder i} and \eqref{eq tilde N}, we have 
    \[
    \frac{\log M_i}{r_i}>h(\mu_i)-\epsilon_i>h(m)-\frac{2\delta}{3}>h(\nu)-\delta,
    \]
    and for each $x\in \bigcup_{j=1}^{M_i}J_{i,j}$ we have 
    \begin{align*}
    &\frac{\log |(f^{r_i})'(x)|}{r_i}<\chi(\mu_i)+\epsilon_i<\chi(m)+\frac{2\delta}{3}< \chi(\nu)+\delta \text{ and }
    \\&\left|\frac{S_{r_i}\phi(x)}{r_i}-\nu(\phi)\right|<
    \left|\frac{S_{r_i}\phi(x)}{r_i}-\mu_i(\phi)\right|+\frac{\delta}{3}+\frac{\delta}{3}<\delta.
    \end{align*}
\end{proof}
For each $N\in\mathbb{N}$ we define 
\[
\mathcal{A}_N:=\bigcup_{i=1}^N\mathcal{J}_i=\bigcup_{i=1}^N\bigcup_{j=1}^{M_i}\{J_{i,j}\}.
\]
Let $N\in\mathbb{N}$.
For each $1\leq i\leq N$ and $J\in \mathcal{J}_{i}$ we set 
\[
r(J):=r_i \text{ and } g_J:=(f^{r(J)}|_{J})^{-1}:B\rightarrow J\subset W_i.
\]
Recall that $\{W_i\}_{i\in\mathbb{N}}$ is pairwise disjoint. Hence, $\mathcal{A}_N$ is pairwise disjoint. Moreover, we have 
\begin{align}\label{eq uniformly expanding F_N}
\inf _{J\in \mathcal{A}_N}\inf_{x\in J} |(f^{r(J)})'(x)|\geq 2 .   
\end{align}
Thus, by an argument similar to that in Remark
\ref{rem coding katok}, the coding map $\Pi_N:\mathcal{A}_N^{\mathbb{N}}\rightarrow \Lambda_N :=\Pi_N(\mathcal{A}_N^{\mathbb{N}})$ defined by 
\[
\{\Pi_N(\tau)\}=\bigcap_{n\in\mathbb{N}} g_{\tau_1}\circ\cdots\circ g_{\tau_n}(B)
\]
is well defined. We endow $\mathcal{A}_N^{\mathbb{N}}$
with the shift metric $d$ defined by 
$
d(\tau,\tau')=\exp(-\inf\{k\ge 1:\,\, \tau_k \neq \tau_k' \})
$ for $\tau,\tau'\in \mathcal{A}_N^{\mathbb{N}}$. Then $\Pi_N$ is a homeomorphism. Moreover, $\Pi_N$ and $\Pi^{-1}_N$ are H\"older continuous. 
We define the induced map $F_N:\Lambda_N\rightarrow \Lambda_N$ by 
\[
F_N|_{\Lambda_N\cap J}:=f^{r(J)} \text{ for each }J\in \mathcal{A}_N.
\]
Then we have 
\[
\Pi_N\circ\sigma_N=F_N\circ\Pi_N. 
\]
Here, $\sigma_N$ denotes the left shift map on $\mathcal{A}_N^{\mathbb{N}}$.
Hence, 
by \eqref{eq uniformly expanding F_N}, the map $F_N$ is uniformly expanding and topologically conjugate to the full shift $\mathcal{A}_N^{\mathbb{N}}$.

For each $1\leq i \leq N$ we define 
\[
D_{i,N}:=\bigcup_{j=1}^{M_i}(J_{i,j}\cap\Lambda_N) \text{ and } Y_N:=\bigcup_{i=1}^N\bigcup_{k=0}^{r_i-1}f^k(D_{i,N}). 
\]
\begin{lemma}\label{lemma uniformly expanding Y_N}
  For each $N\in\mathbb{N}$ the set $Y_N$ is compact, $f(Y_N)=Y_N$ and $Y_N\subset Y_{N+1}$. Moreover, the map $f|_{Y_N}$ is uniformly expanding.   
\end{lemma}

\begin{proof}
Let $N\in\mathbb{N}$.
Since $D_{i,N}$ is compact for each $1\leq i\leq N$, the set $Y_N$ is compact. 
Since $\mathcal{A}_{N}\subset\mathcal{A}_{N+1}$, for each $1\leq i\leq N$ we have $D_{i,N}\subset D_{i,N+1}$. Hence, we obtain $Y_N\subset Y_{N+1}$. 
Moreover, we have 
\[
f(Y_N)=\bigcup_{i=1}^N\bigcup_{k=1}^{r_i-1}f^k(D_{i,N}) \cup \bigcup_{i=1}^N F_N(D_{i,N})=\bigcup_{i=1}^N\bigcup_{k=1}^{r_i-1}f^k(D_{i,N}) \cup \Lambda_N=Y_N.
\]
We will show that $f|_{Y_N}$ is uniformly expanding. By \eqref{eq uniformly expanding F_N}, the compactness of $J\in\mathcal{A}_N$ implies that  
\[
C:=\min \left\{1, \min_{J\in \mathcal{A}_N}\min _{0\leq n\leq m\leq r(J)}\min_{x\in J}|(f^{m-n})'(f^n(x))|\right\}>0,
\]
where $f^0$ denotes the identity map.
Let $y\in Y_N$. Then there exist $1\leq i\leq N$, $0\leq j\leq r_i-1$ and $x_0\in D_{i,N}$ such that $y=f^{j}(x_0)$. For each $\ell\in\mathbb{N}$ we set $x_\ell:=F^\ell_N(x_0)$. For each $\ell\in \mathbb{N}\cup \{0\}$ there exists a unique interval $J_\ell\in \mathcal{A}_N$ such that $x_\ell\in J_\ell$. Since $x_{\ell+1}=F_N(x_{\ell})=f^{r(J_\ell)}(x_\ell)$, for each $\ell\in\mathbb{N}$ we have 
\begin{align*}
x_{\ell}=f^{\sum_{k=0}^{\ell-1}r(J_k)}(x_0).
\end{align*}
We set 
\[
t_0:=-j\text{ and } t_\ell:=-j+\sum_{k=0}^{\ell-1}r(J_k) \text{ for each }\ell\in\mathbb{N}.
\]
Then we have 
\begin{align}\label{eq uniformly expanding 0}
t_0\leq 0< t_1<t_2<\cdots \text{ and }f^{t_\ell}(y)=x_{\ell} \text{ for each }\ell\in\mathbb{N}.
\end{align}
In particular, for each $n\in\mathbb{N}$ we obtain
\begin{align}\label{eq uniformly expanding nLR}
&    n=\sum_{\ell=0}^\infty|[0,n]\cap [t_{\ell},t_{\ell+1}]|\leq L_nR+2R=(L_n+2)R \implies L_n\geq \frac{n}{R}-2,
\\& \text{where } 
R:=\max\{r_i:1\leq  i\leq N\} \text{ and } 
L_n:=\#\{\ell\in\mathbb{N}:0\leq t_\ell< t_{\ell+1}\leq n\}.
\nonumber
\end{align}
Let $n\geq t_1$. There exists a unique number $\ell\in\mathbb{N}$ such that $t_\ell\leq n<t_{\ell+1}$. 
Since $n=t_1+ \sum_{k=1}^{\ell-1}r(J_k)+(n-t_\ell)$, \eqref{eq uniformly expanding 0}, \eqref{eq uniformly expanding F_N} and \eqref{eq uniformly expanding nLR} imply that 
\begin{align*}
    &|(f^n)'(y)|=|(f^{t_1})'(y)|\prod_{k=1}^{\ell-1}|(f^{r(J_k)})'(x_k)||(f^{n-t_\ell})'(x_\ell)|
    \\&\geq C^22^{\ell-1}=C^22^{L_n}\geq  2^{-2}C^2(2^{1/R})^n.
\end{align*}
This completes the proof.
\end{proof}

Let $N\in\mathbb{N}$.
We define $\rho:\Lambda_N\rightarrow \{r_i:1\leq i\leq N\}$ by 
\[
\rho|_{D_{i,N}}\equiv r_i.
\]
Then for each $x\in \Lambda_N$ we have $F_N(x)=f^{\rho(x)}(x)$. For each function $\psi$ on $Y_N$ we define $\tilde \psi: \Lambda_N\rightarrow\mathbb{R}$ by 
\[
\tilde \psi(x):=\sum_{k=0}^{\rho(x)-1} \psi (f^{k}(x)).
\]
For each ergodic measure $\tilde \mu\in M(F_N)$ we define 
\[
\mathcal{L}_N(\tilde \mu):=\frac{1}{\tilde \mu(\rho)}\sum_{i=1}^N\sum_{j=0}^{r_i-1}\tilde \mu|_{D_{i,N}}\circ f^{-j}=\frac{1}{\tilde \mu(\rho)}\sum_{n=0}^\infty\sum_{k=n+1}^\infty\tilde \mu|_{\{\rho=k\}}\circ f^{-n}.
\]
\begin{lemma}\label{lemma lift}
Let $N\in\mathbb{N}$.  We have 
\begin{align}\label{eq surgective L_N}
\mathcal{L}_N(\{\tilde \mu\in M(F_N): \tilde \mu \text{ is ergodic}\})=\{ \mu\in M(f|_{Y_N}):  \mu \text{ is ergodic}\}.
\end{align}
Moreover, for each ergodic measure $\mu\in M(f|_{Y_N})$ with $\mu=\mathcal{L}_N(\tilde \mu)$ for some $\tilde \mu\in M(F_N)$ and each Borel measurable function $\psi: Y_N\rightarrow \mathbb{R}$ with $\mu(|\psi|)<\infty$ we have 
\begin{align}\label{eq Abramov Kac}
   \tilde \mu(\rho) h(\mu)=h(\tilde \mu),\ \tilde \mu(\rho) \mu(\psi)=\tilde \mu(\tilde \psi) \text{ and }
    \tilde \mu(\rho)\chi(\mu)=\tilde \mu(\log |F_N'|).
\end{align}
Furthermore, the map 
\[
\mathcal{L}_N:\{\tilde \mu\in M(F_N): \tilde \mu \text{ is ergodic}\}\rightarrow \{ \mu\in M(f|_{Y_N}):  \mu \text{ is ergodic}\}
\]
is continuous.
\end{lemma}

\begin{proof}
Let $N\in\mathbb{N}$.
Since we have $Y_N\subset \bigcup_{j=0}^{\max\{r_i:1\leq i\leq N\}-1}f^{-j}(\Lambda_N)$, for any $\mu\in M(f|_{Y_N})$ we have $1\leq \max\{r_i:1\leq i\leq N\} \mu(\Lambda_N)$. Thus, for any $\mu\in M(f|_{Y_N})$ we have $\mu(\Lambda_N)>0$. Therefore,  
    \eqref{eq surgective L_N} follows from \cite[Proposition 2.2 and Theorem 2.4]{PesinSenti} (see also \cite{zweimuller}). \eqref{eq Abramov Kac} follows from \cite[Theorem 2.3]{PesinSenti}. The continuity of the map $\mathcal{L}_N$ follows from \eqref{eq Abramov Kac}.
\end{proof}

\begin{rem}\label{rem Y_N property}
    Let $N\in\mathbb{N}$. By \cite[Theorem 4.5]{PesinSenti}, for any H\"older continuous function $\psi: Y_N\rightarrow \mathbb{R}$ there exists a unique equilibrium state $\mu_\psi$ for $\psi$, that is, 
    \[
    \sup\{h(\mu)+\mu(\psi):\mu\in M(f|_{Y_N})\}=h(\mu_\psi)+\mu_\psi(\psi).
    \]
\end{rem}

\begin{lemma}\label{lemma approximation any measure by YN}
Let $\mu\in M(f)$. We assume that $h(\mu)>0$ and $\mu(\phi)\in (\alpha_{\min},\alpha_{\max})$. Then for each $\zeta>0$ there exist $N\in\mathbb{N}$ and $\delta>0$ such that for each $\alpha\in (\mu(\phi)-\delta, \mu(\phi)+\delta)$ there exists an ergodic measure $\nu\in M(f|_{Y_N})$ such that  $\nu(\phi)=\alpha$, $h(\nu)>h(\mu)-\zeta$, $h(\nu)>0$ and $\chi(\nu)<\chi(\mu)+\zeta$.    
\end{lemma}
\begin{proof}
    By Lemma \ref{lemma nonempty B}, there exist $\mu_-, \mu_+\in M(f)$ such that $h(\mu_-)>0$, $h(\mu_+)>0$ and $\mu_-(\phi)<\mu(\phi)<\mu_+(\phi)$. Let 
    \[
    \varepsilon:=\min\{\mu(\phi)-\mu_-(\phi), \mu_+(\phi)-\mu(\phi)\}>0.
    \]
    By Lemma \ref{lemma any measure aproximated for some i}, there exist $i_-,i_+\in\mathbb{N}$ such that $i_-<i_+$, 
    \begin{align}\label{eq average - +}
        &\left|\frac{S_{r_{i_-}}\phi(x)}{r_{i_-}}-\mu_-(\phi)\right|<\frac{\varepsilon}{4} \text{ for all }x\in \bigcup_{j=1}^{M_{i_-}}J_{i_-,j} \text{ and }
        \\&       \left|\frac{S_{r_{i_+}}\phi(x)}{r_{i_+}}-\mu_+(\phi)\right|<\frac{\varepsilon}{4} \text{ for all }x\in \bigcup_{j=1}^{M_{i_+}}J_{i_+,j}.\nonumber
    \end{align}
    We define 
    \[
    \eta:=\min\{\zeta, h(\mu)\}/4.
    \]
    We fix a small number $0<\theta<1$ satisfying 
    \begin{align}\label{eq def theta}
        \theta h(\mu)<\eta \text{ and }\theta C<\eta,
    \end{align}
    where $C:= \log\max\{ 2,  \|f'\|_\infty\}$.
    By Lemma \ref{lemma any measure aproximated for some i},  there exists $N>i_+$ such that for all $x\in \bigcup_{j=1}^{M_N}J_{N,j}$ we have 
    \begin{align}\label{eq def N}
    &\frac{\log M_N}{r_N}>h(\mu)-\eta, \ \frac{\log |(f^{r_N})'(x)|}{r_N}<\chi(\mu)+\eta 
    \\&\text{ and } \left|\frac{S_{r_N}\phi(x)}{r_N}-\mu(\phi)\right|<\frac{\theta \varepsilon}{4}.
  \nonumber  \end{align}
    For each $t\in [0,1]$ we define the probability vector $\{w_i(t)\}_{i=1}^N$ by $w_i(t)=0$ if $i\notin \{N,i_-,i_+\}$,
    \begin{align*}
    w_i(t):=1-\theta \text{ if }i=N,\ w_i(t):=\theta(1-t) \text{ if }i=i_{-} \text{ and } w_i(t):=\theta t \text{ if }i=i_+.    
    \end{align*}
    Let $t\in [0,1]$.
    For each $1\leq i\leq N$ and each $J\in \{J_{i,1},\cdots, J_{i,M_i}\}$ we define 
    \[
    p_i(t):=\frac{w_i(t)/r_i}{\sum_{k=1}^N w_k(t)/r_k} \text{ and }u_J(t):=\frac{p_i(t)}{M_i}.
    \]
    Then, $\{u_J(t)\}_{J\in \mathcal{A}_N}$ is a probability vector. Indeed, we have $u_J(t)\in [0,1]$ for each $J\in \mathcal{A}_N$ and 
    \begin{align*}
        \sum_{J\in \mathcal{A}_N}u_J(t)=\sum_{i=1}^N\sum_{k=1}^{M_i}u_{J_{i,k}}(t)=\frac{\sum_{i=1}^N w_i(t)/r_i}{\sum_{k=1}^N w_k(t)/r_k}=1.
    \end{align*}
 We denote by $\tilde \nu_t'$ the Bernoulli measure on $\mathcal{A}_N^\mathbb{N}$ associated with the probability vector $\{u_J(t)\}_{J\in \mathcal{A}_N}$. Define
 \[
 \tilde \nu_t:=\tilde \nu_t'\circ \Pi_N^{-1} \text{ and }\nu_t:=\mathcal{L}_N(\tilde \nu_t).
 \]
 Notice that 
 \begin{align*}
     \tilde \nu_t(\rho)=\sum_{i=1}^N\sum_{k=1}^{M_i} r_i u_{J_{i,k}}(t)=\sum_{i=1}^N p_i(t) r_i=\frac{1}{\sum_{k=1}^Nw_k(t)/r_k}.
 \end{align*}
 Therefore, for each $1\leq i\leq N$ we have
 \begin{align}\label{eq nu t return}
     \frac{p_i(t) r_i}{\tilde \nu_t(\rho)}=\frac{w_i(t)}{\tilde \nu_t(\rho) \sum_{k=1}^Nw_k(t)/r_k}=w_i(t).
 \end{align}
 By Lemma \ref{lemma lift}, we have $\nu_t\in M(f|_{Y_N})$ and  $\nu_t$ is ergodic. 
 By \eqref{eq Abramov Kac}, \eqref{eq nu t return}, \eqref{eq def N} and \eqref{eq def theta}, we have 
 \begin{align}\label{eq nu t entropy}
     & h(\nu_t)=\frac{-\sum_{i=1}^N\sum_{k=1}^{M_i}u_{J_{i,k}}(t)\log u_{J_{i,k}}(t)}{\tilde \nu_t(\rho)}
     \\&\nonumber
=\frac{-\sum_{i=1}^Np_i(t)\log p_i(t)}{\tilde \nu_t(\rho)}+\sum_{i=1}^N w_i(t)\frac{\log M_i}{r_i}\geq \sum_{i=1}^N w_i(t)\frac{\log M_i}{r_i}
\\&\geq \nonumber (1-\theta)\frac{\log M_N}{r_N}>h(\mu)-\eta-\theta h(\mu)+\theta \eta>h(\mu)-2\eta>\max\{0,h(\mu)-\zeta\}.
 \end{align}
By \eqref{eq Abramov Kac}, \eqref{eq nu t return}, Ruelle's inequality \eqref{eq ruelle's ine}, \eqref{eq def N} and \eqref{eq def theta} we have 
\begin{align}\label{eq nu t lyapunov}
&\chi(\nu_t)=\frac{\sum_{i=1}^N\sum_{k=1}^{M_i}r_i\int_{J_{i,k}} {\log |(f^{r_i})'|}/{r_i}d\tilde \nu_{t}}{\int \rho d\tilde \nu_t}
\\&\nonumber\leq\frac{r_Np_N(t)(\chi(\mu)+\eta) }{\int \rho d\tilde \nu_t}+\frac{r_{i_-}p_{i_-}(t)C}{\int \rho d\tilde \nu_t}+\frac{r_{i_+}p_{i_+}(t)C}{\int \rho d\tilde \nu_t}
=w_N(t)(\chi(\mu)+\eta)\\&+C(w_{i_-}(t)
+w_{i_+}(t))
= (1-\theta)(\chi(\mu)+\eta)+\theta C\leq \chi(\mu)+2\eta<\chi(\mu)+\zeta.\nonumber
\end{align}
    By Lemma \ref{lemma lift}, the function $A:[0,1]\rightarrow \mathbb{R}$ defined by 
 \[
 A(t):=\nu_t(\phi)
 \]
 is continuous. By \eqref{eq Abramov Kac}, \eqref{eq average - +} and \eqref{eq def N}, for all $t\in [0,1]$ we have 
 \begin{align*}
     &\left|A(t)-(w_{i_-}(t)\mu_-(\phi)+w_{i_+}(t)\mu_+(\phi)+w_N(t)\mu(\phi))\right|=
     \\&\left|\frac{\sum_{i=1}^N\int _{D_{i,N}}S_{r_{i}}\phi d\tilde \nu_t-\left(p_{i_-}(t)r_{i_-}\mu_{-}(\phi)+p_{i_+}(t)r_{i_+}\mu_{+}(\phi)+p_{N}(t)r_{N}\mu(\phi)\right)}{\int \rho d\tilde \nu_t}\right|
     \\&\leq\Bigg|
     \frac{
     \int_{D_{i_-,N}}
     (S_{r_{i_-}}\phi-r_{i_-}\mu_{-}(\phi))d\tilde \nu_t}{\int \rho d\tilde\nu_t}\Bigg|
     +\Bigg|
     \frac{
     \int_{D_{i_+,N}}(S_{r_{i_+}}\phi-r_{i_+}\mu_{+}(\phi))d\tilde \nu_t}{\int \rho d\tilde\nu_t}\Bigg|
     \\&+\left|\frac{\int_{D_{N,N}}(S_{r_N}\phi-r_N\mu(\phi))d\tilde\nu_t}{p_{N}(t)r_{N}/w_{N}(t)}\right|
     <\frac{\varepsilon}{4}(\theta(1-t)+\theta t+(1-\theta)\theta)<\frac{\theta \varepsilon}{2}.
 \end{align*}
    In particular, we have 
    \begin{align*}
        &A(0)<\theta\mu_-(\phi)+(1-\theta)\mu(\phi)+\frac{\theta \varepsilon}{2}\leq \mu(\phi)-\theta\varepsilon+\frac{\theta \varepsilon}{2}=\mu(\phi)-\frac{\theta\varepsilon}{2} \text{ and }
        \\&A(1)>\theta\mu_{+}(\phi)+(1-\theta)\mu(\phi)-\frac{\theta \varepsilon}{2}\geq \mu(\phi)+\theta \varepsilon-\frac{\theta\varepsilon}{2}=\mu(\phi)+\frac{\theta \varepsilon}{2}.
    \end{align*}
    Hence, we obtain $A(0)<\mu(\phi)-
    \theta\varepsilon/2<\mu(\phi)+\theta\varepsilon/2<A(1)$. 
    Let $\delta:=\theta\varepsilon/2$.
    Since $A$ is continuous, for any $\alpha\in (\mu(\phi)-\delta,\mu(\phi)+\delta)$ there exists $t\in [0,1]$ such that $\nu_t(\phi)=\alpha$. By \eqref{eq nu t entropy} and \eqref{eq nu t lyapunov}, we obtain the desired result.
    \end{proof}

\subsection{Proof of Theorem \ref{thm uniform compact appro spectrum}}\label{sec proof of uniform}
Throughout this section, we use the same notation
as in the previous section.

\begin{prop}
    For each $\alpha\in (\alpha_{\min},\alpha_{\max})$ we have 
    \[
    D(\alpha)=\sup_{N\in\mathbb{N}}\sup \left\{\frac{h(\mu)}{\chi(\mu)}:\mu\in M(f|_{Y_N})\text{ is ergodic and }\mu(\phi)=\alpha\right\}.
    \]
\end{prop}
\begin{proof}
    Let $\alpha\in (\alpha_{\min},\alpha_{\max})$. By Lemma \ref{lemma uniformly expanding Y_N}, we have 
    \begin{align}\label{eq lowerbound apporximation}
    D(\alpha)\geq \sup_{N\in\mathbb{N}}\sup \left\{\frac{h(\mu)}{\chi(\mu)}:\mu\in M(f|_{Y_N})\text{ is ergodic and }\mu(\phi)=\alpha\right\}.
    \end{align}
    By Lemma \ref{lemma nonempty B}, we have 
    \[
    D(\alpha)=\sup\left\{\frac{h(\mu)}{\chi(\mu)}:\mu\in M(f),\ h(\mu)>0,\ \mu(\phi)=\alpha\right\}.
    \]
    Hence, by Lemma \ref{lemma approximation any measure by YN}, we obtain  
    \[
    D(\alpha)\leq \sup_{N\in\mathbb{N}} \sup\left\{\frac{h(\mu)}{\chi(\mu)}:\mu\in M(f|_{Y_N})\text{ is ergodic and }\mu(\phi)=\alpha\right\}.
    \]
    Combining this with \eqref{eq lowerbound apporximation}, we obtain the desired result.
\end{proof}

For any $\nu\in M(f)$ we define 
\[
\dim_H(\nu):=\inf\{\dim_H(X): \text{$X$ is a Borel subset of $I$ satisfying }\nu(X)=1\}.
\]
Recall that the function $f'$ is H\"older continuous. Thus, by \cite[Theorem 1]{Hofbauer}, we have the following: 
For each ergodic measure $\nu\in M(f)$ with $\chi(\nu)>0$ we have 
\begin{align}\label{eq Young formula}
    \dim_H(\nu)=\frac{h(\nu)}{\chi(\nu)}.
\end{align}

\emph{Proof of Theorem \ref{thm uniform compact appro spectrum}.}
We first show \eqref{eq main uniform weak}.
Let $K\subset \hyperbolic\cap (\alpha_{\min},\alpha_{\max})$ be a nonempty compact set. Let $\epsilon>0$ and let $\alpha_0\in K$. By \eqref{eq conditional variational principle strong} in  Theorem \ref{thm main}, we have $D(\alpha_0)=b(\alpha_0)$. Thus, by Lemma \ref{lemma nonempty B}, there exists $\mu\in M(f)$ such that $\mu(\phi)=\alpha_0$, $h(\mu)>0$, and $h(\mu)/\chi(\mu)>b(\alpha_0)-\epsilon/3$. We take a small number $0<\zeta<h(\mu)$ satisfying 
\begin{align}\label{eq appro}
\frac{h(\mu)-\zeta}{\chi(\mu)+\zeta}>b(\alpha_0)-\frac{2\epsilon}{3}.
\end{align}
By Lemma \ref{lemma approximation any measure by YN}, there exist $N(\alpha_0)\in\mathbb{N}$ and $\delta(\alpha_0)>0$ such that for each $\alpha\in (\alpha_0-\delta(\alpha_0), \alpha_0+\delta(\alpha_0))$ there exists an ergodic measure $\nu_\alpha\in M(f|_{Y_{N(\alpha_0)}})$ such that $\nu_\alpha(\phi)=\alpha$, $h(\nu_\alpha)>h(\mu)-\zeta$ and  $\chi(\nu_\alpha)<\chi(\mu)+\zeta$. By Birkhoff's ergodic theorem, \eqref{eq Young formula} and \eqref{eq appro}, for each $\alpha\in (\alpha_0-\delta(\alpha_0), \alpha_0+\delta(\alpha_0))$ we obtain
\begin{align}\label{eq lower bound by measure}
\dim_H(Y_{N(\alpha_0)}\cap B(\alpha))\geq \frac{h(\nu_\alpha)}{\chi(\nu_\alpha)}>\frac{h(\mu)-\zeta}{\chi(\mu)+\zeta}>b(\alpha_0)-\frac{2\epsilon}{3}.
\end{align}
By Corollary \ref{cor main}, $b$ is continuous on $\hyperbolic\cap (\alpha_{\min},\alpha_{\max})$. Thus, there exists $0<\eta(\alpha_0)<\delta(\alpha_0)$ such that for all $\alpha\in (\alpha_0-\eta(\alpha_0), \alpha_0+\eta(\alpha_0))\cap K$ we have $|b(\alpha_0)-b(\alpha)|<\epsilon /3$. Hence, by \eqref{eq lower bound by measure}, for each $\alpha\in (\alpha_0-\eta(\alpha_0), \alpha_0+\eta(\alpha_0))\cap K$ we have 
\begin{align}\label{eq appro main}
0\leq b(\alpha)-\dim_H(Y_{N(\alpha_0)}\cap B(\alpha))<\epsilon.
\end{align}
Since $K$ is compact and $K= \bigcup_{\alpha_0\in K}(\alpha_0-\eta(\alpha_0), \alpha_0+\eta(\alpha_0))\cap K$, there exist $\ell\in \mathbb{N}$ and $\{\alpha_1,\cdots, \alpha_\ell\}\subset K$ such that 
\[
K\subset \bigcup_{j=1}^\ell(\alpha_j-\eta(\alpha_j), \alpha_j+\eta(\alpha_j)).
\]
For each $\alpha\in K$ we fix $j(\alpha)$ such that $\alpha\in (\alpha_{j(\alpha)}-\eta(\alpha_{j(\alpha)}), \alpha_{j(\alpha)}+\eta(\alpha_{j(\alpha)}))$.
Define $N_0:=\max\{N(\alpha_j):1\leq j\leq \ell\}$.
Since $Y_N\subset Y_{N+1}$ for each $N\in\mathbb{N}$, \eqref{eq appro main} implies that for each $N\geq N_0$ and $\alpha\in K$
\begin{align*}
    0\leq b(\alpha)-\dim_H(Y_N\cap B(\alpha))\leq b(\alpha)-\dim_H(Y_{N(\alpha_{j(\alpha)})}\cap B(\alpha))<\epsilon.
\end{align*}
This completes the proof of \eqref{eq main uniform weak}.

Next, we show \eqref{eq main uniform strong}. We assume that $H=1$. Then by Theorem \ref{thm main} and Corollary \ref{cor main}, the function $b$ is continuous on $(\alpha_{\min},\alpha_{\max})\setminus \partial\bad$ and for each $\alpha\in (\alpha_{\min},\alpha_{\max})\setminus \partial\bad$ we have $b(\alpha)=D(\alpha)$. Hence, by an argument similar to that used to prove
\eqref{eq main uniform weak}, we obtain
\eqref{eq main uniform strong}.
\qed

\subsection*{Acknowledgments}
This work was supported by the JSPS KAKENHI 25KJ1382.

\bibliographystyle{abbrv}
\bibliography{reference}

\end{document}